\documentclass[11pt]{article}
\usepackage{amscd,amssymb,stmaryrd}
\usepackage{amsmath}
\usepackage{graphicx}
\usepackage{subcaption}
\usepackage[utf8]{inputenc}
\usepackage[export]{adjustbox}
\usepackage{wrapfig}
\usepackage{amstext}
\usepackage{pstricks,pst-node,pst-plot,pst-coil}
\usepackage{amsthm}
\usepackage{amsmath}
\usepackage{color}
\usepackage{mathtools}
\usepackage{hyperref}
\usepackage[english]{babel}
\usepackage{booktabs}
\usepackage{mathrsfs}
\usepackage{comment}
\usepackage{float}
\usepackage[linesnumbered,ruled,vlined, ruled,vlined]{algorithm2e}

\theoremstyle{definition}

\theoremstyle{remark}

\def\eg{e.g., }
\def\ie{i.e., }

\title{Accelerated replica exchange stochastic gradient Langevin diffusion enhanced Bayesian DeepONet for solving noisy parametric PDEs}
\author{Guang Lin\footnote{Department of Mathematics and Mechanical Engineering, Purdue University, West Lafayette, IN 47906, USA}, ~ Christian Moya\footnote{Department of Mathematics, Purdue University, West Lafayette 47906, IN, USA}, ~Zecheng Zhang \footnote{Department of Mathematics, Purdue University, West Lafayette 47906, IN, USA}}
\begin{document}
\maketitle

\begin{abstract}
    The Deep Operator Networks~(DeepONet) is a fundamentally different class of neural networks that we train to approximate nonlinear operators, including the solution operator of parametric partial differential equations~(PDE). DeepONets have shown remarkable approximation and generalization capabilities even when trained with relatively small datasets. However, the performance of DeepONets deteriorates when the training data is polluted with noise, a scenario that occurs very often in practice. To enable DeepONets training with noisy data, we propose using the Bayesian framework of replica-exchange Langevin diffusion. Such a framework uses two particles, one for exploring and another for exploiting the loss function landscape of DeepONets. We show that the proposed framework's exploration and exploitation capabilities enable (1) improved training convergence for DeepONets in noisy scenarios and (2) attaching an uncertainty estimate for the predicted solutions of parametric PDEs. In addition, we show that replica-exchange Langeving Diffusion (remarkably) also improves the DeepONet's mean prediction accuracy in noisy scenarios compared with vanilla DeepONets trained with state-of-the-art gradient-based optimization algorithms (\eg Adam). To reduce the potentially high computational cost of replica, in this work, we propose an accelerated training framework for replica-exchange Langevin diffusion that exploits the neural network architecture of DeepONets to reduce its computational cost up to 25\% without compromising the proposed framework's performance. Finally, we illustrate the effectiveness of the proposed Bayesian framework using a series of experiments on four parametric PDE problems. 
\end{abstract}

\section{Introduction}
Many problems in Science and Engineering require solving parametric PDEs \cite{bhattacharya2020model, khoo2021solving, fresca2021comprehensive, zhang2020learning}. That is, they require the repeated evaluation of a PDE model for a given distribution of inputs. For example, in petroleum engineering applications, engineers seek to calculate, using the porous media equations, the pressure field of the oil based on its permeability. In practice, the value of the oil permeability varies frequently; hence, one needs to calculate the oil pressure field for a distribution of fast-varying permeabilities~\cite{zhang2020learning,chung2020multi, chung2021multi}. The traditional numerical frameworks require intense computations to solve these parametric PDEs. Furthermore, the computational cost for these traditional frameworks increases even more when the problem is time-dependent and multiscale~\cite{chung2021computational, chetverushkin2021computational, efendiev2021hei}. 

To reduce the computational cost for solving parametric PDEs, many works have proposed using deep neural networks, which transfer most of the cost to offline training.  However, the data required to train these networks is often polluted with noise or numerical errors in practice. In these noisy scenarios, deep neural networks often have poor generalization performance because the gradient-based algorithms used to train the networks tend to chase the noise. As a result, the trained model will fail if used in practice.  

This work focuses on Deep Operator Networks~(DeepONet)~\cite{lu2021learning}, a fundamentally different class of neural networks that learns to approximate nonlinear operators, i.e., mappings between infinite-dimensional function spaces. DeepONets are constructed to realize and exploit the power of the universal approximation theorem of operators~\cite{chen1995universal}. To that end, a DeepONet uses two sub-networks: one  is called the Branch net, which encodes the input function at a fixed number of locations; and another one is called the Trunk net, which encodes the locations where one evaluates the output function. Such an architecture enables DeepONets to learn nonlinear operators efficiently from relatively small datasets. The remarkable capabilities of DeepONets have been reported in many application areas, including electroconvection~\cite{cai2021deepm}, chemistry~\cite{lin2021operator}, economics~\cite{leite2021DeepONets}, and for solving parametric PDEs~\cite{lu2021learning}. These previous works, however, have only trained DeepONets using data without noise. In real applications, data is often polluted with noise. So, a framework that can handle noisy data during DeepONets' training must be developed.

Many works have studied the problem of dealing with noisy datasets and avoiding overfitting. Among these works, Bayesian frameworks~\cite{welling2011bayesian} that use the Langevin diffusion~\cite{welling2011bayesian, chen2014stochastic} have provided promising results. For example, in~\cite{welling2011bayesian}, the authors introduced the Stochastic Gradient Langevin diffusion algorithm to train neural networks. They demonstrated that the samples generated by their proposed algorithm converge to the target posterior distribution. Replica-exchange MCMC and (its stochastic gradient variants) were proposed~\cite{lin2021multi, chen2020accelerating, deng2020non} to accelerate the Langevin diffusion. In replica-exchange frameworks, instead of using one particle to sample from the posterior, one uses two particles with different temperature parameters. The high-temperature particle explores the landscape of the loss function, which enables escaping from local minima, while the low temperature exploits the same landscape, enforcing local convergence~\cite{hwang1980laplace, zhang2017hitting}. The two particles can be exchanged according to a swapping probability. If this swapping probability is designed so that the Markov process is reversible, then the low-temperature particle may attain the global optima faster. 

Replica-exchange frameworks provide us with excellent convergence performance. They, however, double the computational on each training iteration. Hence, one needs to balance having fast convergence (\ie fewer iteration steps) and haven high per-iteration computational cost. To tackle such a problem, in our previous work~\cite{lin2021multi}, we proposed the multi-variance stochastic gradient Langevin diffusion~(m-reSGLD) framework, which assumes that the estimators for the two particles have different accuracy. Hence, one can use solvers of different accuracy when using replica-exchange frameworks. 

In this work, our objective is to develop a Bayesian framework that can train DeepONets to approximate the solution operator of parametric PDEs using noisy measurements. To this end, we develop a replica-exchange stochastic gradient Langevin diffusion~(reSGLD) algorithm tailored for training DeepONets. Furthermore, we design an accelerated training regime for replica-exchange frameworks. Such a training regime results in a multi-variance replica-exchange SGLD (m-reSGLD) algorithm that reduces the cost of reSGLD. We summarize the contributions of this work next. 

\begin{itemize}
    \item We design a Bayesian framework to approximate the solution operator of parametric PDEs using DeepONets in scenarios where the training data is polluted with noise. To the best of the authors' knowledge, this is the first paper that addresses such a problem. In the proposed framework, the Bayesian DeepONet represents the prior for the trainable parameters, while replica-exchange Stochastic Gradient Langevin Diffusion~(reSGLD) enables estimating the posterior, which we use to estimate the uncertainty of DeepONet predictions. 
    \item We demonstrate that the replica-exchange algorithm ability to escape local minima, using one particle to explore and another particle to exploit the DeepONet's loss function landscape, enables the proposed framework to provide improved training convergence for DeepONets in noisy scenarios when compared to vanilla DeepONets trained with gradient-based optimization algorithms (\eg Adam). Remarkably, the proposed framework also provides a more accurate mean prediction performance than vanilla DeepONets. 
    \item We also propose a novel accelerated Bayesian training regime that reduces the computational cost of the replica-exchange algorithm up to 25\% without compromising the DeepONet's predictive performance. More precisely, for the particle that explores the loss landscape, we randomly choose to train one of the DeepONet's sub-networks (the Branch Net or the Trunk Net) while keeping the other fixed.  
    \item Finally, we test the effectiveness of the proposed Bayesian framework using four parametric PDE problems. In particular, we compare the training convergence speed, mean prediction accuracy, and uncertainty quantification performance of reSGLD, the proposed accelerated reSGLD (denoted as m-reSGLD), and Adam with dropout (a non-Bayesian framework for uncertainty quantification).
\end{itemize}
The rest of the paper is organized as follows. In Section~\ref{sec:background}, we provide a brief review of Deep Operator Network~(DeepONet). Section~\ref{sec:bayesian-deeponet} introduces the Bayesian setting of our problem and details the replica-exchange Stochastic Gradient Langevin diffusion~(reSGLD) algorithm. We present the proposed accelerated replica-exchange algorithm in Section~\ref{sec:accelerated-training}. In Section~\ref{sec:numerical}, we test the performance of the proposed Bayesian framework using a series of experiments on four parametric PDE problems. Finally, Section~\ref{sec:conclusion} concludes this work.

\section{Background Information} \label{sec:background}
In this work, we propose a Bayesian data-driven deep learning framework to approximate the solution operator of the prototypical parametric PDE
\begin{align}  \label{eq:PDE}
(\mathcal{L}_a s)(x) = f(x), \qquad x \in D, 
\end{align}
where $D \subset \mathbb{R}^d$ is a bounded open set, \textcolor{black}{$a:D \to \mathbb{R}$ is a parameter entering the definition of the nonlinear operator~$\mathcal{L}_a$}, and $s: D \to \mathbb{R}$ is the solution to the PDE (given appropriate boundary conditions). To approximate the solution operator, we use the deep operator neural network~(DeepONet) framework introduced in~\cite{lu2021learning}, which we review next.

\subsection{Review of DeepONet} \label{subsec:review-deeponet}
\begin{figure}[t!]
\centering
\includegraphics[width=.55\textwidth, height=3.0cm]{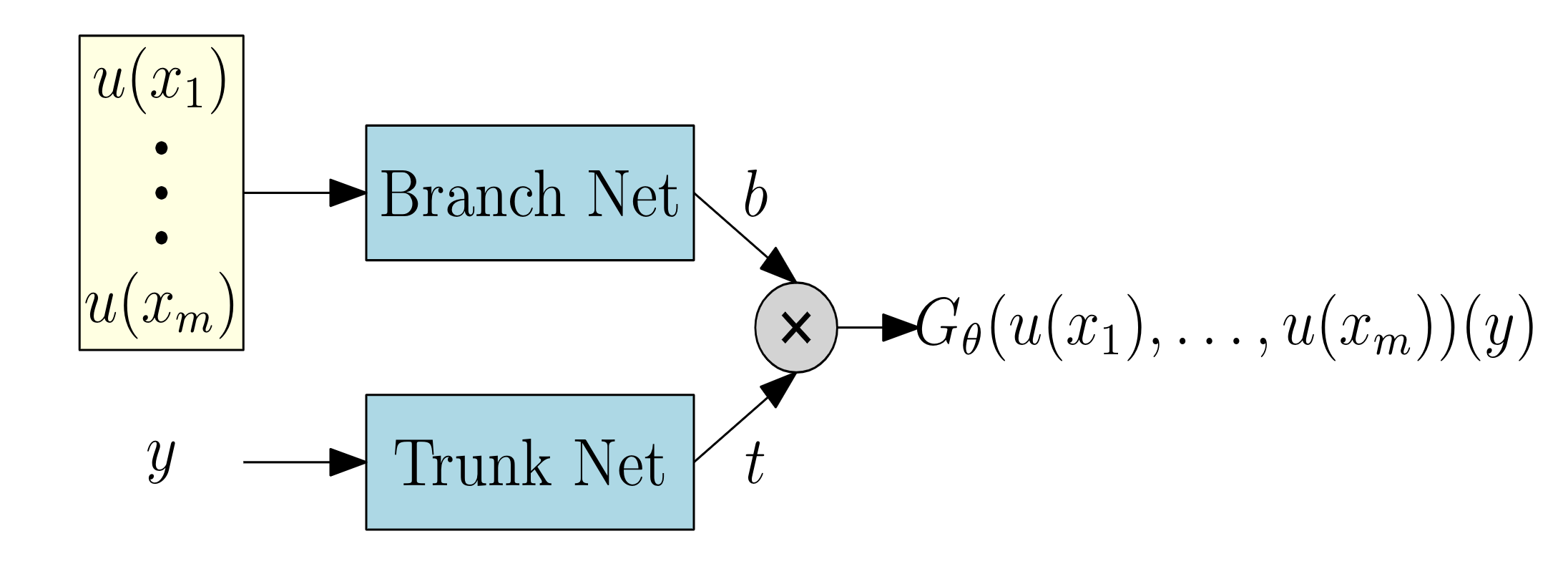}
\caption{The DeepONet architecture. The crossed node on the right indicates the DeepONet output, which we obtain by taking the inner product between the output features of the Branch~(b) and Trunk~(t) Nets.}
\label{fig:deepONet}
\end{figure}

Let $G^\dagger$ denote the nonlinear operator arising as the solution operator of the parametric PDE~\eqref{eq:PDE}. This operator~$G^\dagger$ maps an input function~$u$ (corresponding, for example, to the parameter $a$, the forcing~$f$, or initial/boundary conditions) to an output function~$G^\dagger(u)$ (corresponding to the solution $s$ of the parametric PDE). Let $y \in Y$ denote a point in the output function domain~$Y \subset \mathbb{R}^d$. Then, the goal of the DeepONet~$G_\theta$, with trainable parameters~$\theta \in \mathbb{R}^p$, is to approximate the operator~$G^\dagger(u)(y)$ at $y \in Y$. To this end, the DeepONet~$G_\theta$ uses the architecture depicted in Figure~\ref{fig:deepONet}, consisting of two sub-networks referred to as the Branch Net and Trunk Net. 

The \textit{Branch} Net processes the input function information. Let $(x_1,\ldots, x_m)$ denote points in the domain of $u$ (we refer to these points as the \textit{sensors}), such that $(u(x_1),\ldots,u(x_m))$ is a discrete representation of the input $u$. The Branch Net takes this discretized $u$ as the input and outputs a vector of features~$b \in \mathbb{R}^q$. On the other hand, the \textit{Trunk} Net processes points in the domain~$Y$ of the output function. To this end, the Trunk Net takes $y \in  Y$ as the input and outputs a vector of features~$t \in \mathbb{R}^q$. Note that since the Trunk Net's output~$t$ solely depends on the input coordinates~$y$, it is natural to interpret the components of~$t$ as a collection of basis functions defined on~$Y$, \ie
$$
t =(\varphi_1(y),\ldots,\varphi_q(y)).
$$
The output of the DeepONet then combines the output features from the Branch Net~$b$ and the Trunk Net~$t$ using an inner product:
\begin{align} \label{eq:deeponet}
    G_\theta \left(u(x_1),\ldots,u(x_m)\right)(y) := \langle b,t \rangle = \sum_{i=1}^{q} b_i \cdot \varphi_i(y).
\end{align}
From the above, one can interpret the output of the Branch Net~$b$ as the trainable coefficients for the basis functions~$t$ produced by the Trunk Net. To simplify our notation, in the rest of this work, we omit writing the DeepONet explicit dependency on the discretized input and use the simplified notation~$G_\theta(u)(y)$.

Finally, we train the DeepONet~$G_\theta$ to approximate the nonlinear solution operator~$G^\dagger$ by minimizing a mean square loss function on the training dataset $\left\{u_i,y_i,G^\dagger(u_i)(y_i)\right\}_{i=1}^N$, \ie
\begin{align} \label{eq:deeponet-loss}
    \mathcal{L}(\theta) = \frac{1}{N} \sum_{i = 1}^N \left|G_\theta(u_i)(y_i) - G^\dagger(u_i)(y_i)\right|^2,
\end{align}
In this work, our aim is to desing a Bayesian training framework for DeepONets that enables us to (1) better approximate~$G^\dagger$ and (2) estimate the solution/prediction uncertainty. We refer to this framework as the \textit{Bayesian DeepONet}.

\section{The Bayesian DeepONet} \label{sec:bayesian-deeponet}
To derive the proposed Bayesian DeepONet framework, we consider the scenario when our available operator training targets correspond to scattered noisy measurements of $G^\dagger(u)(y)$, \ie \textcolor{black}{$\{\tilde{G}^\dagger_i\}_{i=1}^N$}. In particular, we assume these measurements are independently Gaussian distributed centered at the latent true operator target value, \ie 
$$
\tilde{G}^\dagger_i = G^\dagger(u_i)(y_i) + \epsilon_i, \qquad i=1,2,\ldots,N.
$$
Here, $\epsilon_i$ is an independent Gaussian noise with zero mean and \textit{known} standard deviation~$\sigma$, \ie $\epsilon_i \sim \mathcal{N}(0, \sigma^2)$. Let $\mathcal{D}:=\{\tilde{G}^{\dagger}_i\}_{i=1}^N$ denote the noisy dataset of targets; we can then calculate the likelihood as:
\begin{align} \label{eq:likelihood}
    P(\mathcal{D}|\theta) = \prod_{i=1}^N P(\tilde{G}^\dagger_i|\theta) = \prod_{i=1}^N \frac{1}{\sqrt{2 \pi \sigma^2}} \exp \left( - \frac{(G_\theta(u_i)(y_i) - \tilde{G}^\dagger_i)^2}{2 \sigma^2} \right).
\end{align}
To compute the posterior distribution, we use Bayes' theorem, \ie 
$$
P(\theta | \mathcal{D}) = \frac{P(\mathcal{D}| \theta) P(\theta)}{P(\mathcal{D})} \propto P(\mathcal{D}| \theta) P(\theta),
$$
where ``$\propto$'' denotes equality up to a constant. In practice, computing $P(\mathcal{D})$ is usually analytically intractable. Hence, to sample from the posterior~$P(\theta | \mathcal{D})$, we usually use the unnormalized expression $P(\mathcal{D}| \theta) P(\theta)$.

To predict and estimate the uncertainty of a solution trajectory, denoted as~$\{G^\dagger(u)(y):y \in Y\}$, for a given input~$u$ and over the mesh of points~$y \in Y$, we must sample from the posterior~$P(\theta | \mathcal{D})$ and obtain an $M$-ensemble of DeepONet parameters, \ie~$\{\theta_k\}_{k=1}^{M}$. Then, using this ensemble, we can obtain statistics from simulated trajectories $\{G_{\theta_k}(u)(y):y \in Y\}_{k=1}^M$. In this work, we compute the mean and standard deviation of~$\{G_{\theta_k}(u)(y):y \in Y\}_{k=1}^M$. We use the mean to predict the true solution $\{G^\dagger(u)(y):y \in Y\}$ and the standard deviation to quantify the predictive uncertainty. 

To sample from the posterior and obtain the $M$-ensemble $\{\theta_k\}_{k=1}^{M}$, in this work, we use replica exchange stochastic Langevin diffusion~(re-SGLD)~\cite{deng2020non}, which we review in the next section and detail in Algorithm~\ref{alg:re-SGLD}.
\begin{algorithm}[t]
\DontPrintSemicolon
\SetAlgoLined
\textbf{Require:} initial DeepONet parameters~$\theta^1_0, \theta^2_0$, learning rates~$\eta^1_k, \eta^2_k$, temperatures~$\tau_1,\tau_2$, stochastic gradient variances~$\sigma_1^2,\sigma_2^2$, and two constants $a_1,a_2 > 0$ satisfying $a_1+a_2 = 1$.\\
\For{$k = 1,\ldots,N$}{
  \textbf{sampling step:}\;
  sample $B_k^1 \sim \mathcal{N}(0, I)$ and $B_k^2 \sim \mathcal{N}(0, I)$\; 
  $\qquad \theta^1_{k+1}=\theta^1_k - \nabla \hat{U}(\theta^1_k) + \sqrt{2 \tau_1 \eta_k^1} B_k^1$\;
  $\qquad \theta^2_{k+1}=\theta^2_k - \nabla \hat{U}(\theta^2_k) + \sqrt{2 \tau_2 \eta_k^2} B_k^2$\;
  \textbf{swapping step:}\;
  Generate a uniform random number~$u \in [0,1]$\;
  let $\tau_\delta = \frac{1}{\tau_1} - \frac{1}{\tau_2}$ and compute\;
  $\qquad \hat{r} = e^{\tau_\delta \left((a_1-a_2)(\hat{U}_1(\theta^1_{k+1}) - \hat{U}_2(\theta^2_{k+1})) -(a_1 \sigma_1 + a_2 \sigma_2)^2 \tau_\delta \right)}$ \;
\If{$u < \hat{r}$}{
  swap $\theta^1_{k+1}$ and $\theta^2_{k+2}$ \;
}
  }
  Calculate $\{G_{\theta^1_{N+1-k}}(u)(y):y \in Y\}_{k=1}^{M}$ as prediction samples of the true posterior trajectory $\{G^\dagger(u)(y): y \in Y\}$.\;
 \caption{Replica Exchange Stochastic Gradient Langevin Diffusion}
 \label{alg:re-SGLD}
\end{algorithm}

\subsection{Review of the re-SGLD} \label{sec:re-SGLD}
To derive re-SGLD (Algorithm~\ref{alg:re-SGLD}), we start from the Langevin diffusion, which reads:
\begin{align}  \label{eq:langevin}
    d \theta_t = -\nabla U(\theta_t)+\sqrt{2\tau}dB_t,
\end{align}
where $\theta_t \in \mathbb{R}^p$ is called the particle and corresponds to the target trainable parameters of the DeepONet, $B_t$ is the Brownian motion in $\mathbb{R}^p$, $\tau > 0$ is the temperature parameter, and $U: \mathbb{R}^p \to \mathbb{R}$ is an energy function defined as:  
\begin{align} \label{eq:energy}
    U(\theta) \sim -\log P(\theta) - \sum_{i = 1}^N \log P(\tilde{G}^\dagger_i|\theta),
\end{align}
where $P(\theta)$ denotes the prior and $P(\tilde{G}^\dagger_i|\theta)$ the likelihood~\eqref{eq:likelihood}.

To obtain the target parameters, we must solve an optimization problem. More precisely, we must find $\theta^* \in \mathbb{R}^p$ that minimizes the energy~\eqref{eq:energy}. To this end, one can show that, under proper assumptions for the energy function~$U$, the sequence of targets~$\{\theta_t\}_{t \ge 0}$ converges to the target distribution $\pi_{\mathcal{D}} \sim \exp \left( - \frac{U(\theta)}{\tau} \right)$ as $t \to \infty$~\cite{chiang1987diffusion}. Furthermore, when the temperature parameter~$\tau$ is small enough, $\pi_\mathcal{D} (\theta)$ distributes around the minimizer~$\theta^*$~\cite{hwang1980laplace}. Hence, one can obtain the minimizer~$\theta^*$ of the energy~\eqref{eq:energy} by sampling from the target distribution~$\pi_\mathcal{D}(\theta)$.

In addition, when the temperature~$\tau$ is low, the particle tends to reach stationary points with fast local convergence~\cite{zhang2017hitting}, by \textit{exploiting} the landscape of the energy function~$U$. On the other hand, when the temperature~$\tau$ is high, the particle tends to traverse globally, exploring the landscape of the energy function~$U$. This exploitation-exploration dilemma has led researchers to design algorithms (\eg~re-SGLD) that achieve superior performance by balancing exploration and exploitation. 

In this work, we use exploration and exploitation by considering two Langevin equations with separated particles $\theta^1_t$ and $\theta^2_t$ and different temperatures $\tau_1$ (low) and $\tau_2$ (high), \ie 
\begin{subequations} \label{eq:re-LD}
\begin{align}
    d \theta_t^1 = -\nabla U(\theta_t^1)+\sqrt{2\tau_1}dB_t^1,\\
    d \theta_t^2 = -\nabla U(\theta_t^2)+\sqrt{2\tau_2}dB_t^2.
\end{align}
\end{subequations}
Then, by swapping the particles $\theta^1_t$ and $\theta^2_t$, with controlled rate~$r(\theta^1_t, \theta^2_t)$, one can help the low-temperature particle~$\theta^1_t$ to escape the local minima and achieve improved and accelerated convergence. 

The theoretical support for the previous claims were established in the convergence result of~\cite{chen2020accelerating}. That is, the \textit{replica exchange Langevin diffusion}~(re-LD) algorithm, which simulates the dynamics~\eqref{eq:re-LD} and allows the swapping of particles, \ie $(\theta^1_{t+dt},\theta^2_{t+dt}) = (\theta^2_{t+dt},\theta^1_{t+dt})$, with rate:
\begin{align} \label{eq:rate}
    r(\theta^1_t, \theta^2_t) = e^{\tau_{\delta}(U(\theta_t^1) - U(\theta_t^2) ) },
\end{align}
converges to the invariant distribution with density:
\begin{align} \label{eq:target-distr-re-LD}
    \pi(\theta^1, \theta^2) \propto e^{-\frac{U(\theta^1)}{\tau_1}-\frac{U(\theta^2)}{\tau_2}}.
\end{align}
\subsection{Errors in the energy function} \label{sub-sec:errors-energy}
In practice, computing the energy function~$U$ may fail due to, for example, errors in the dataset~$\mathcal{D}$ or an insufficient number of samples of~$\mathcal{D}$ to represent the likelihood and cover the target output space. Moreover, when the size~$N$ of the dataset is large, it is expensive to compute~$U$. Hence, the authors in~\cite{deng2020non} proposed to approximate the energy function~$U$ using a mini-batch of data $\{\tilde{G}^\dagger_{s_i}\}_{i=1}^n \subset \mathcal{D}$. In any case, an unbiased estimator of the energy function is then 
\begin{align} \label{eq:estimator}
    \hat{U}(\theta) = -\log P(\theta) - \frac{N}{n}\sum_{i = 1}^n \log P(\tilde{G}^\dagger_{s_i}|\theta).
\end{align}
Moreover, we assume the estimator~\eqref{eq:estimator} follows the normal distribution
$$
\hat{U}(\theta) \sim \mathcal{N}(U(\theta), \sigma_e^2),
$$
where $\sigma^2_e$ is the variance of the estimator. To simulate from the dynamics~\eqref{eq:langevin}, we discretize them and obtain
\begin{subequations} \label{eq:discrete-langevin}
\begin{align} 
     \theta_{k+1}^1 = \theta_{k}^1-\nabla \hat{U}(\theta_{k}^1)+\sqrt{2\tau_1\eta_k}B_k^1,\\
     \theta_{k+1}^2 = \theta_{k}^2-\nabla \hat{U}(\theta_{k}^2)+\sqrt{2\tau_2\eta_k}B_k^2,
\end{align}
\end{subequations}
where $\eta_k^1, \eta_k^2$ are the positive learning rates and $\hat{U}_1(\theta^1)\sim\mathcal{N}(U(\theta^1), \sigma_1^2), \hat{U}_2(\theta^2)\sim\mathcal{N}(U(\theta^2), \sigma_1^2)$ are the energy function estimators of the two particles. However, using the unbiased estimators for the energy, $\hat{U}(\theta^1)$ and $\hat{U}(\theta^2)$, in re-LD with discretized dynamics~\eqref{eq:discrete-langevin} leads to a large bias for the estimator of the swapping rate~$r(\theta^1_t, \theta^2_2)$ defined in~\eqref{eq:rate}. 

To remove the bias from the swaps, we allow the particles swapping $(\theta_{k+1}^1, \theta_{k+1}^2 ) = (\theta_{k+1}^2, \theta_{k+1}^1 ) $ with the following unbiased rate estimator~\cite{deng2020non}, 
\begin{align}
      \Hat{r} = e^{\tau_\delta\bigg(a_1\big(\Hat{U}_1(\theta^1)-\Hat{U}_1(\theta^2)\big)+a_2\big(\Hat{U}_2(\theta^1)-\hat{U}_2(\theta^2) \big)
    - (a_1\sigma_1  + a_2\sigma_2)^2\tau_\delta
     \bigg) }
\end{align}
where $a_1+a_2 = 1$ are two positive constants. The replica exchange algorithm then converges to the target invariant distribution with this new unbiased rate estimator~$\hat{r}$. We refer to this algorithm as the \textit{replica exchange stochastic gradient Langevin diffusion}~(reSGLD) if $\sigma_1 = \sigma_2$ (see Algorithm~\ref{alg:re-SGLD}); otherwise, we refer to the algorithm as the \textit{multi-variance replica exchange stochastic gradient Langevin diffusion}~(m-reSGLD).

\section{Accelerated Bayesian training of DeepONets} \label{sec:accelerated-training}
In this section, we propose exploiting the DeepONet's architecture to develop an accelerated Bayesian training strategy that (1) reduces the time reSGLD (Algorithm~\ref{alg:re-SGLD}) uses to train the DeepONet parameters and (2) achieves a performance comparable to reSGLD. Such training strategy results in a multi-variance replica-exchange SGLD~(m-reSGLD) algorithm (Algorithm~\ref{alg:m-re-SGLD}), which we describe next.

\begin{algorithm}[t]
\DontPrintSemicolon
\SetAlgoLined
\textbf{Require:} initial DeepONet parameters~$\theta^1_0, \theta^2_0$, learning rates~$\eta^1_k, \eta^2_k$, temperatures~$\tau_1,\tau_2$, stochastic gradient variances~$\sigma_1^2,\sigma_2^2$, constants $a_1,a_2 > 0$ satisfying $a_1+a_2 = 1$, and training control parameter~$c \in [0,1]$.\\
\For{$k = 1,\ldots,N$}{
  \textbf{sampling step:}\;
  sample $B_k^1 \sim \mathcal{N}(0, I)$ and $B_k^2 \sim \mathcal{N}(0, I)$\; 
  $\qquad \theta^1_{k+1}=\theta^1_k - \nabla \hat{U}(\theta^1_k) + \sqrt{2 \tau_1 \eta_k^1} B_k^1$\;
  Generate a uniform random number $\gamma \in [0,1]$\;
  \eIf{$\gamma < c$}{
  $\qquad \theta^{2,\text{BN}}_{k+1}=\theta^{2,\text{BN}}_k - \nabla \hat{U}(\theta^{2,\text{BN}}_k) + \sqrt{2 \tau_2 \eta_k^2} B_k^2$\;
  $\qquad \theta^{2,\text{TN}}_{k+1}=\theta^{2,\text{TN}}_k$\:
  }
  {
  $\qquad \theta^{2,\text{BN}}_{k+1}=\theta^{2,\text{BN}}_k$\;
  $\qquad \theta^{2,\text{TN}}_{k+1}=\theta^{2,\text{TN}}_k - \nabla \hat{U}(\theta^{2,\text{TN}}_k) + \sqrt{2 \tau_2 \eta_k^2} B_k^2$\;
  }
  \textbf{swapping step:}\;
  Generate a uniform random number~$u \in [0,1]$\;
  let $\tau_\delta = \frac{1}{\tau_1} - \frac{1}{\tau_2}$ and compute\;
  $\qquad \hat{r} = e^{\tau_\delta \left((a_1-a_2)(\hat{U}_1(s^1_{k+1}) - \hat{U}_2(s^2_{k+1})) -(a_1 \sigma_1 + a_2 \sigma_2)^2 \tau_\delta \right)}$ \;
\If{$u < \hat{r}$}{
  swap $\theta^1_{k+1}$ and $\theta^2_{k+2}$ \;
}
  }
  Calculate $\{G_{\theta^1_{N+1-k}}(u)(y):y \in Y\}_{k=1}^{M}$ as prediction samples of the true posterior trajectory $\{G^\dagger(u)(y): y \in Y\}$.\;
 \caption{Accelerated Replica Exchange Stochastic Gradient Langevin Diffusion}
 \label{alg:m-re-SGLD}
\end{algorithm}
Many theoretical results have showed that over-parametrized neural networks, trained with gradient descent-based methods, converge to zero-training loss exponentially fast with the network parameters hardly changing~\cite{du2018gradient, du2019gradient}. Motivated by this convergence result, other works have studied whether the same holds for under-parametrized neural networks. Their results~\cite{chizat2018lazy} show that exponential convergence to zero-training loss can happen, with parameters hardly varying, but depends on an implicit scaling. Motivated by these results, in this work, we propose a more efficient Bayesian framework that trains a less parametrized DeepONet, \ie a DeepONet that keeps some of its parameters hardly varying during training.    

Intuitively, our proposed framework follows from the fact that empirically (see~\cite{lu2021learning}), DeepONets converge exponentially fast to zero-training loss. Hence, we expect that a DeepONet that keeps, during training, some of its parameters hardly varying should also converge to zero-training loss provided an appropriate training regime. In this paper, such a regime corresponds to training the Branch Net or the Trunk Net randomly. We refer to this training regime as the \textit{accelerated Bayesian training of DeepONets}, which we use to reduce the computational cost of reSGLD.

The reSGLD has demonstrated excellent convergence capabilities when used for training neural networks. However, compared to vanilla SGLD~\cite{welling2011bayesian}, the computational cost of reSGLD doubles. To reduce the computational cost of reSGLD, we adopt the idea proposed in~\cite{lin2021multi} and develop m-reSGLD (Algorithm~\ref{alg:m-re-SGLD}); that is, an algorithm that uses (1) a full training regime for the low temperature particle~$\theta^1$, which \textit{exploits} the landscape of the energy function~$U$, and (2) a accelerated training regime for the high temperature particle~$\theta^2$, which \textit{explores} the landscape of~$U$. To this end, we first fully train both particles for a fixed number of burn-in epochs. After burn-in, we let the high temperature particle enter the accelerated training regime. Let us split the high temperature particle into its Branch Net and Trunk Net components, \ie $\theta^2 = \{\theta^{2,\text{BN}},\theta^{2,\text{TN}}\}$. Then, to control which sub-network parameters ($\theta^{2,\text{BN}}$ or $\theta^{2,\text{TN}}$) are trained more, we define the training control parameter $c \in [0,1]$. \textcolor{black}{If $c > 0.5$, we train more the Branch Net parameters~$\theta^{2,\text{BN}}$. If, on the other hand, $c < 0.5$}, we train more the Trunk Net paramters~$\theta^{2,\text{TN}}$. 

We select the parameters that we train more based on the DeepONet architecture. As explained in Section~\ref{subsec:review-deeponet}, the DeepONet architecture, depicted in Figure~\ref{fig:deepONet} and whose output is given by~\eqref{eq:deeponet}, enables us to interpret the Branch Net output~$b$ as trainable coefficients for the basis functions~$t$ produced by the Trunk Net. Hence, DeepONets, as proved in~\cite{kovachki2021neural}, fall within the class of linear approximation methods because they approximate the nonlinear space of the operator $G^\dagger$ via the linear space generated by the Trunk Net, \ie $\text{span}\{\varphi_1(y), \dots, \varphi_q(y)\}$. Thus, in practice, we choose to train more the Branch Net because its output represents the coefficients of the linear space generated by the Trunk Net.

Let us conclude this section with the following result. If the networks are identical, then the accelerated training regime reduces at most half of the computational cost for $\theta^2$. Thus, overall, m-reSGLD reduces at most $25 \%$ of the reSGLD computational cost. Finally, one should note that, due to the accelerated regime, the variances of the particles are no longer the same. As result, we use the name multi-variance reSGLD for Algorithm~\ref{alg:m-re-SGLD}.

\section{Numerical experiments}
\label{sec:numerical}
In this section, we test the accuracy and efficiency of the proposed Bayesian training framework using four examples: an anti-derivative operator, a gravity pendulum system, a diffusion reaction system, and an advection diffusion system. To train the DeepONets that will approximate the solution operator of each example, we use three frameworks. First,  we use the state-of-the-art gradient descent optimization framework Adam~(Adam)~\cite{kingma2014adam}. We then compare the performance of Adam with the proposed replica exchange stochastic gradient Langevin diffusion~(reSGLD) and the accelerated replica exchange framework~(m-reSGLD). During training, each framework will be tested systematically using datasets of noisy output targets~$\{\tilde{G}_i^\dagger\}_{i=1}^N$ with different noise values. 

\textit{Metrics.} To test the performance of the frameworks, we compute the $L^1$ and $L^2$ relative errors of a test trajectory. That is, for a given input $(u(x_1),..., u(x_m))$, we predict using DeepONets the solution of this test trajectory at a collection of selected mesh points $y \in Y_m \subset Y$. Denote the mean predicted solution of the test trajectory as $v_u := \{G_\theta(u)(y):y \in Y_m\}$ and the true solution as $w_u :=\{G^\dagger(u)(y):y \in Y_m\}$, then we calculate the relative errors as follows:
\begin{align}
    e_1 = \frac{\|v_u-w_u\|_1}{\|w_u\|_1}100, \quad e_2 = \frac{\|v_u-w_u\|_2}{\|w_u\|_2}100. 
\end{align}
In our experiments, we compute mean value of the relative errors for 100 test trajectories selected from outside the training dataset.

We also verify how the three frameworks handle the noisy targets and estimate the predictive uncertainty. To this end, for every test trajectory, we construct a $95\%$ ($2 \sigma$) confidence interval. For reSGLD and m-reSGLD the confidence interval is constructed using the $M$ prediction samples $\{G_{\theta^1_{N+1-k}}(u)(y): y\in Y_m\}_{k=1}^M$ obtained, respectively, from Algorithm~\ref{alg:re-SGLD} and Algorithm~\ref{alg:m-re-SGLD}. To enable Adam to estimate uncertainty, we adopt the non-Bayesian framework of dropout~\cite{gal2016dropout} and use the strategy described in~\cite{zhang2019quantifying}. To measure how well these confidence intervals quantify uncertainty, we compute the ratio of the true trajectory that is within the confidence interval, \ie
\begin{align}
    e_3 = \frac{\text{\# of points of the predicted solution~$w_u$ in the confidence interval}}{\text{\# of mesh points of the true solution~$w_u$}} 100.
\end{align}
\subsection{Experiment 1}
In this experiment, we use the proposed Bayesian framework to train a DeepONet that approximates the solution operator of the following ordinary differential equation~(ode):
\begin{align} \label{eq:ode}
    \frac{ds}{dt} = u(t), \quad t\in[0, 1],
\end{align}
with the initial condition is $u(0) = 0$. Note that the solution operator for~\eqref{eq:ode} corresponds to the anti-derivative operator:
$$
s(t) = \int_{0}^t u(\tau) d\tau, \quad t \in [0,1].
$$

\begin{figure}[h!]
\centering
\includegraphics[scale = 0.5]{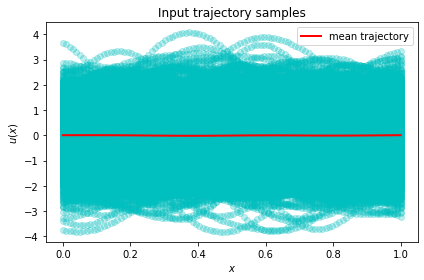}
\caption{Distribution of the input samples~$u$ simulated from the GRF $\mathcal{G}(0, k_l(x_1, x_2))$ with RBF kernel and length-scale $l= 0.2$.} 
\label{input_distribution}
\end{figure}

Here \textcolor{black}{$u(t)$ denotes an external force that drives the response of the ode system~\eqref{eq:ode}. We sample this external force from the following mean-zero Gaussian Random Field~(GRF):}
$$
u \sim \mathcal{G}(0, k_l(x_1, x_2)),
$$
where the covariance kernel $k_l(x_1, x_2) = \exp(-||x_1, x_2||^2/2l^2)$ is the radial-basis function~(RBF) kernel with length-scale parameter $l>0$. In this experiment, we set the length-scale parameter to $l= 0.2$. Figure~\ref{input_distribution} illustrates the distribution of input samples simulated from the GRF and used to train the DeepONets. More specifically, the Branch Net takes the discretized version~$(u(x_1),\ldots,u(x_m))$ of these input samples~$u$ during training. Here, we discretized $u$ using $m=100$ sensors.  

\begin{figure}[h!]
\centering
\includegraphics[scale = 0.5]{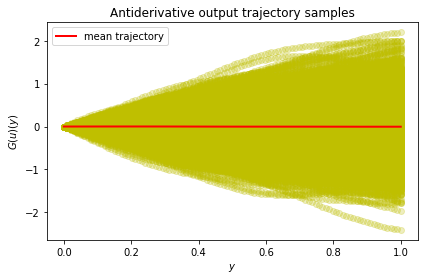}
\caption{Distribution of the operator targets~$G^\dagger(u)(y)$ obtained from the distribution of inputs depicted in Figure~\ref{input_distribution} and used to train the anti-derivative example.} 
\label{anti_output_distribution}
\end{figure}

In addition, we present in Figure~\ref{anti_output_distribution} the distribution of the true operator targets associated with the input samples~$u$. We can observe that these targets deviate largely from the mean. Thus, any predictive framework must capture the operator response over the whole output target space. Such a requirement makes this problem very challenging.

\subsubsection{Small noise}
\begin{figure}[h!]
\centering
\includegraphics[scale = 0.4]{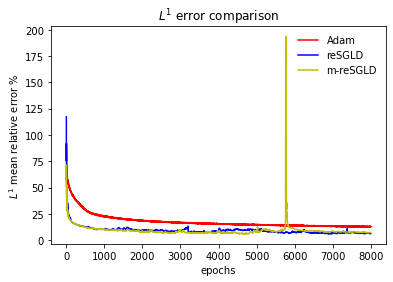}
\includegraphics[scale = 0.4]{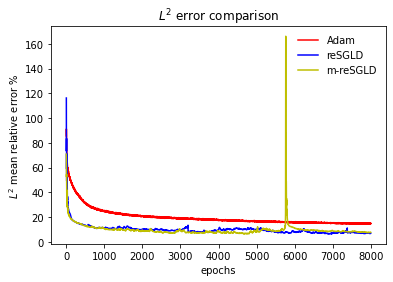}
\includegraphics[scale = 0.4]{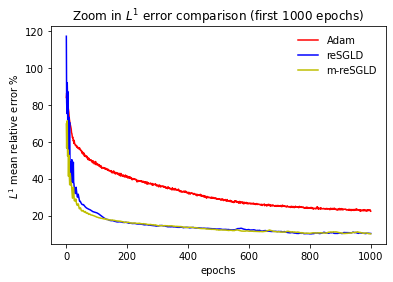}
\includegraphics[scale = 0.4]{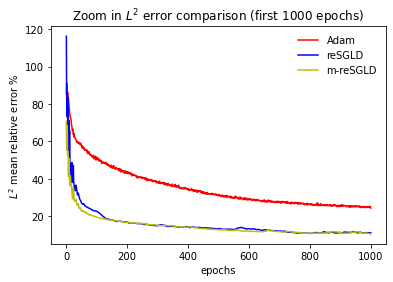}
\includegraphics[scale = 0.4]{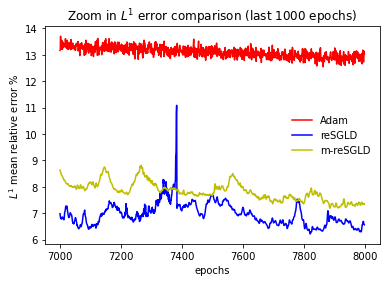}
\includegraphics[scale = 0.4]{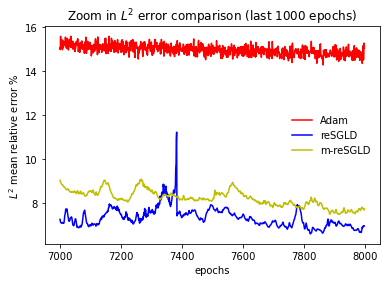}
\caption{Convergence results for the anti-derivative example with noise $\mathcal{N}(0, 0.01^2)$. The $e_1$ errors are presented on the left while $e_2$ errors are shown on the right. Top row: entire history. Middle row: zoom in the first $1,000$ epochs. Last row: zoom in the last $1,000$ epochs. The average errors after the burn in are presented in the Table~\ref{anti_error_small_table}.} 
\label{anti_error_small}
\end{figure}

We consider first training the DeepONet using operator targets~$\{\tilde{G}_1^\dagger\}$ with small noise. That is, we assume the noise~$\epsilon_i$ follows the normal distribution $\mathcal{N}(0, 0.01^2)$.
We train the DeepONet for $8,000$ epochs using the three frameworks. And, at each epoch, we compute the mean of the relative errors~$e_1$ and $e_2$ of 100 test trajectories. Figure~\ref{anti_error_small} shows how these errors converge during training. Furthermore, Table~\ref{anti_error_small_table} presents the average mean relative errors the burn-in epochs.

\begin{table}[h!]
\centering
\begin{tabular}{||c c c||} 
\hline
frameworks & $e_1$ & $e_2$ \\ [0.5ex] 
\hline
Adam & $13.099151$ & $14.961829$ \\ [0.5ex]
\hline
reSGLD  & $6.953939$ & $7.35423$  \\ [0.5ex]
\hline
m-reSGLD   & $7.82799$ & $8.234294$ \\ [0.5ex]
\hline
\end{tabular}
\caption{ Anti-derivative example with noise $\mathcal{N}(0, 0.01^2)$. The mean errors after burn-in for all three frameworks.}
\label{anti_error_small_table}
\end{table}

From Figure~\ref{anti_error_small} and Table~\ref{anti_error_small_table}, we can observe that the proposed reSGLD and m-reSGLD converge much faster than Adam. Moreover, the average relative errors after burn-in show that the mean predictive performance of the proposed Bayesian framework outperforms the state-of-the-art Adam. Such improved accuracy illustrates that reSGLD and m-reSGLD handle much better noisy training data.

\begin{figure}[h!]
\centering
\includegraphics[scale = 0.3]{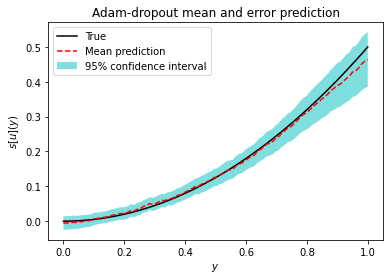}
\includegraphics[scale = 0.3]{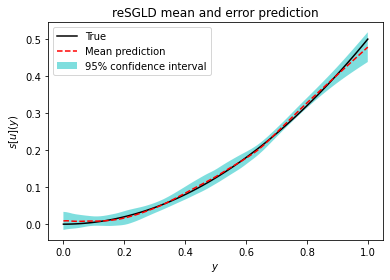}
\includegraphics[scale = 0.3]{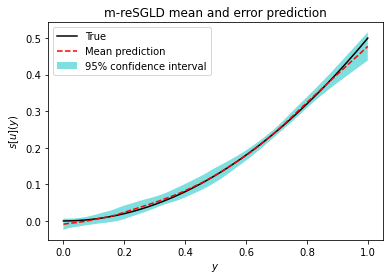}
\caption{Confidence intervals for the anti-derivative example with noise $\mathcal{N}(0, 0.01^2)$. Left: Adam. Middle: reSGLD. Right: m-reSGLD.  In this example, $e_3$ for all three frameworks are $100\%$.} 
\label{anti_confidence_small}
\end{figure}

Besides fast training convergence and improved mean predictive performance, we have to check if the proposed frameworks can capture the predictive uncertainty of the problem. To this end, we illustrate in Figure~\ref{anti_confidence_small}, for a randomly selected test trajectory, the confidence interval and the error~$e_3$ for Adam with dropout, reSGLD, and m-reSGLD. We observe that the three frameworks capture the true trajectory within their confidence intervals. Adam, however, seems to overestimate its predictive uncertainty, producing a wider confidence band. 

\begin{figure}[h!]
\centering
\includegraphics[scale = 0.5]{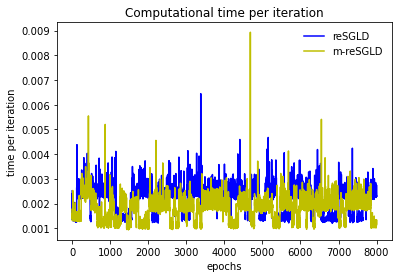}
\caption{Anti-derivative example with noise $\mathcal{N}(0, 0.01^2)$. Computer running time per iteration for reSGLD and m-reSGLD. The average computation time of m-reSGLD is around $80.1\%$ of reSGLD framework.} 
\label{anti_time_small}
\end{figure}

We also verify how much of the computational cost of reSGLD is reduced when using m-reSGLD. To this end, during training, we record the per-iteration computational time for reSGLD and m-reSGLD (see Figure~\ref{anti_time_small}). The results show that the computational time for m-reSGLD corresponds to $80.1\%$ the computational time of reSGLD. This reduction of approximate $20\%$ of the cost of reSGLD is close to the theoretical maximum of~$25\%$. Furthermore, m-reSGLD does not deteriorate the predictive performance of the DeepONet (see Table~\ref{anti_error_small_table}).

\subsubsection{Increased noise}
\begin{figure}[h!]
\centering
\includegraphics[scale = 0.4]{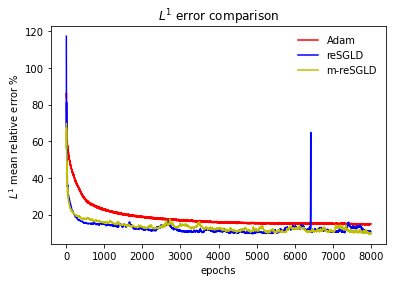}
\includegraphics[scale = 0.4]{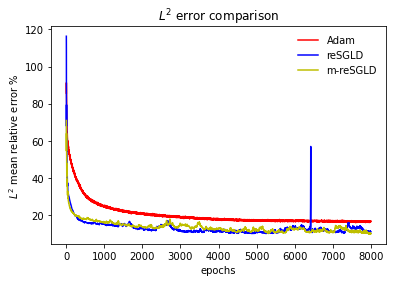}
\includegraphics[scale = 0.4]{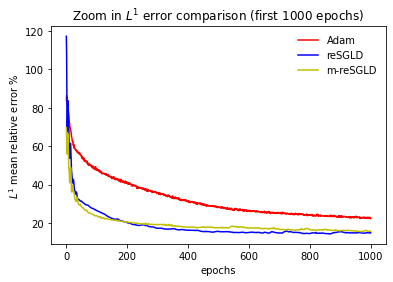}
\includegraphics[scale = 0.4]{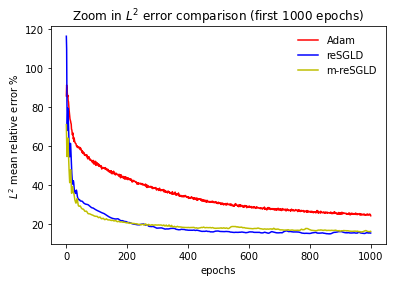}
\includegraphics[scale = 0.4]{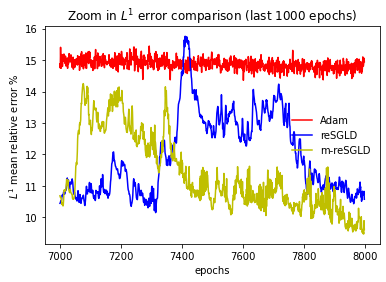}
\includegraphics[scale = 0.4]{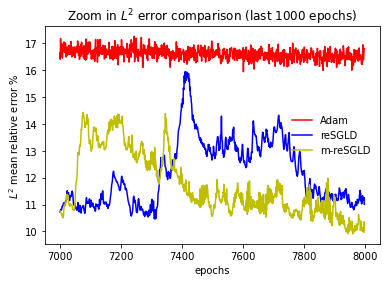}
\caption{Convergence results for the anti-derivative with noise $\mathcal{N}(0, 0.05^2)$. The $e_1$ errors are presented on the left while $e_2$ errors are shown on the right. Top row: entire history. Middle row: zoom in the first $1,000$ epochs. Last row: zoom in the last $1,000$ epochs. The average errors after the burn in are presented in the Table~\ref{anti_error_big_table}.} 
\label{anti_error_big}
\end{figure}

We now verify the performance of the proposed frameworks when the noise associated to the operator output targets~$\{\tilde{G}_i^\dagger\}$ increases. In particular, we assume the increased noise follows the normal distribution $\epsilon_i \sim \mathcal{N}(0, 0.05^2)$. Figure~\ref{anti_error_big} depicts the training convergence results for the the relative errors~$e_1$ and $e_2$, and Table~\ref{anti_error_big_table} presents the average errors.

\begin{table}[h!]
\centering
\begin{tabular}{||c c c||} 
\hline
frameworks & $e_1$ & $e_2$ \\ [0.5ex] 
\hline
Adam & $14.877717$ & $16.611967$ \\ [0.5ex]
\hline
reSGLD  & $11.970458$ & $12.182682$  \\ [0.5ex]
\hline
m-reSGLD   & $11.367707$ & $11.708988$ \\ [0.5ex]
\hline
\end{tabular}
\caption{Anti-derivative example  with noise $\mathcal{N}(0, 0.05^2)$. The mean errors after burn-in for all three frameworks.}
\label{anti_error_big_table}
\end{table}

We observe from Figure~\ref{anti_error_big} and Table~\ref{anti_error_big_table} that the convergence performance of the three frameworks deteriorates in this increased noise scenario. Our proposed frameworks, however, provide a faster training convergence and have much better mean predictive capability than Adam. 

\begin{figure}[h!]
\centering
\includegraphics[scale = 0.35]{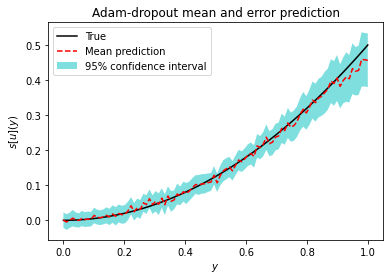}
\includegraphics[scale = 0.35]{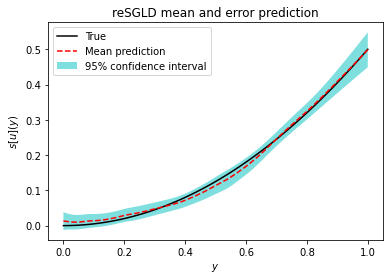}
\includegraphics[scale = 0.35]{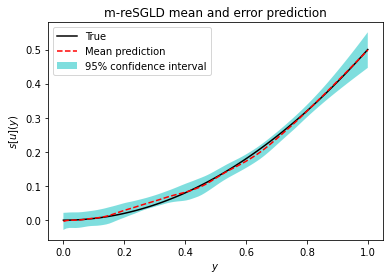}
\caption{Confidence intervals for the anti-derivative example with noise $\mathcal{N}(0, 0.05^2)$. Left: Adam. Middle: reSGLD. Right: m-reSGLD.  In this example, $e_3$ of Adam, reSGLD and m-reSGLD are all $100\%$.} 
\label{anti_confidence_big}
\end{figure}

We also construct the confidence intervals for the three frameworks. Figure~\ref{anti_confidence_big} shows that, as expected, the predictive uncertainty (\ie the width of the confidence band) estimated by the proposed framework has increased in this scenario. The confidence band for Adam, however, has becomes more noisy and uncertain. Also, we note that the three frameworks capture the true solution trajectory, \ie they have~$e_3=100\%$. Thus, we conclude that the proposed frameworks are more effective in predicting uncertainty in this increased noise scenario.

\begin{figure}[h!]
\centering
\includegraphics[scale = 0.5]{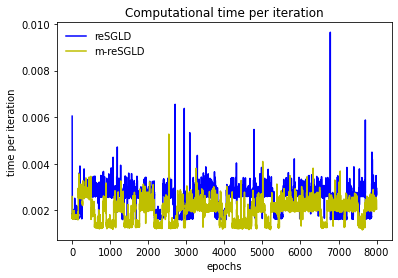}
\caption{Anti-derivative example with noise $\mathcal{N}(0, 0.05^2)$. Computer running time per iteration for reSGLD and m-reSGLD. The average computation time of m-reSGLD is around $76.9\%$ of reSGLD framework.} 
\label{anti_time_big}
\end{figure}

We conclude this experiment by illustrating (see Figure~\ref{anti_time_big}) how much computational cost is saved when using m-reSGLD instead of reSGLD. The average per-iteration time (\ie the computational cost) for m-reSGLD represents $76.9\%$ of the per-teration time used by reSGLD. Moreover, for this increased noise scenario, m-reSGLD is not only more efficient but also provides a better mean predictive performance (see Table~\ref{anti_error_big_table}).

\subsection{Experiment 2}
\begin{figure}[h!]
\centering
\includegraphics[scale = 0.5]{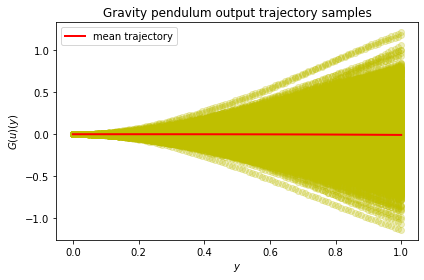}
\caption{Distribution of the operator targets~$G^\dagger(u)(y)$ obtained from the distribution of inputs depicted in Figure~\ref{input_distribution} and used to train the gravity pendulum.} 
\label{pendulum_output_distribution}
\end{figure}

In this experiment, we use the proposed frameworks to train a DeepONet that approximates the solution operator for the following gravity pendulum with external force,
\begin{subequations}
\begin{align}
    \frac{ds_1}{dt} &= s_2,~t \in [0,1], \\
    \frac{ds_2}{dt} &= -k \sin(s_1) + u(t),
\end{align}
\end{subequations}
with an initial condition~$(s_1(0), s_2(0)) = (0,0)$, and $k=1$. We simulate the discretized external force~$u(t)$ using a mean-zero GRF with RBF kernel (length-scale parameter $l=0.2$) and $m=100$ sensors. We illustrate in Figure~\ref{pendulum_output_distribution} the distribution of the true operator targets associated with the input samples~$u$ and used to train the gravity pendulum system.

\subsubsection{Small noise}
\begin{figure}[h!]
\centering
\includegraphics[scale = 0.4]{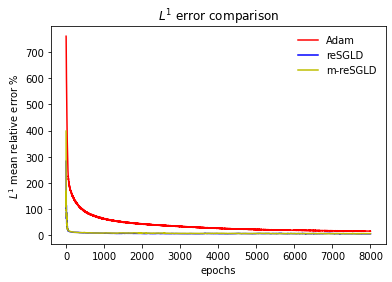}
\includegraphics[scale = 0.4]{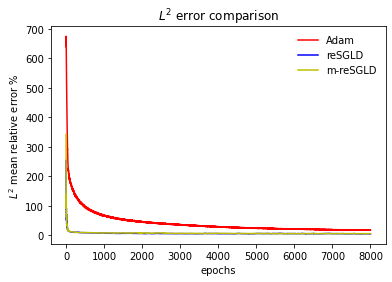}
\includegraphics[scale = 0.4]{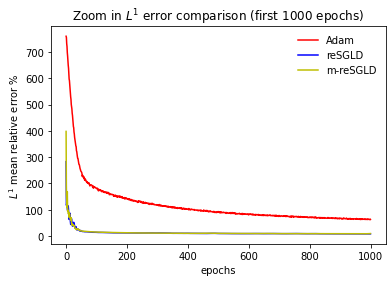}
\includegraphics[scale = 0.4]{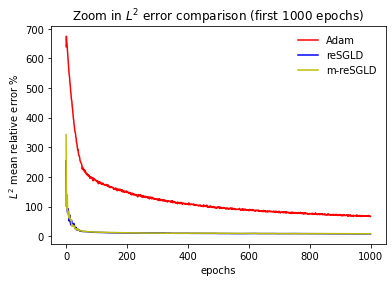}
\includegraphics[scale = 0.4]{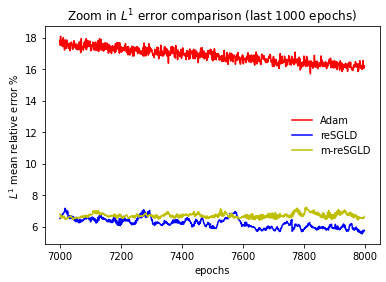}
\includegraphics[scale = 0.4]{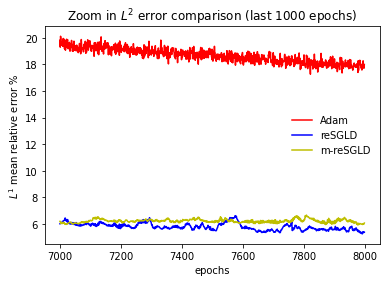}
\caption{Convergence results for the gravity-pendulum with noise $\mathcal{N}(0, 0.01^2)$. The $e_1$ errors are shown on the left while $e_2$ errors are shown on the right. Top row: entire history. Middle row: zoom in the first $1,000$ epochs. Last row: zoom in the last $1,000$ epochs. The average errors after burn-in are presented in Table~\ref{pendulum_error_small_table}. } 
\label{pendulum_error_small}
\end{figure}

We start by verifying the performance of the proposed frameworks using a small noise scenario. In particular, let us assume the noise~$\epsilon_i$ for the target outputs follows the normal distribution~$\epsilon_i \sim \mathcal{N}(0,0.01^2)$. We train the DeepONet for $8,000$ epochs using the three frameworks. Figure~\ref{pendulum_error_small} shows the convergence results for the errors~$e_1$ and $e_2$ during training. And Table~\ref{pendulum_error_small_table} presents the average of these errors after the burn-in epochs. The results show that the proposed frameworks converge much faster during training and provide better mean prediction performance than the state-of-the-art optimization algorithm Adam.

\begin{table}[h!]
\centering
\begin{tabular}{||c c c||} 
\hline
frameworks & $e_1$ & $e_2$ \\ [0.5ex] 
\hline
Adam & $16.8571$ & $18.6582$ \\ [0.5ex]
\hline
reSGLD  & $6.2147$ & $5.7965$  \\ [0.5ex]
\hline
m-reSGLD   & $6.6997$ & $6.2231$ \\ [0.5ex]
\hline
\end{tabular}
\caption{Gravity-pendulum example with noise $\mathcal{N}(0, 0.01^2)$. The mean errors after burn-in for all three frameworks.}
\label{pendulum_error_small_table}
\end{table}

\begin{figure}[h!]
\centering
\includegraphics[scale = 0.35]{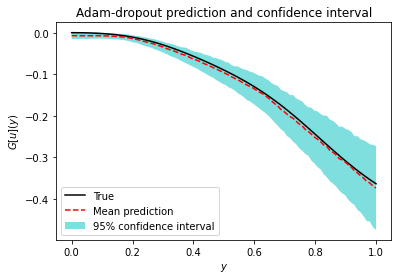}
\includegraphics[scale = 0.35]{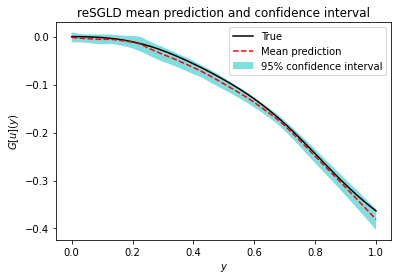}
\includegraphics[scale = 0.35]{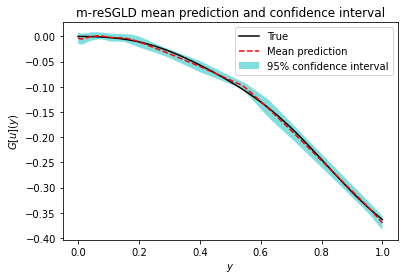}
\caption{Confidence intervals for the gravity pendulum example with noise $\mathcal{N}(0, 0.01^2)$. Left: Adam. Middle: reSGLD. Right: m-reSGLD.  In this example, the error $e_3$ for Adam is $71\%$, for reSGLD is $100\%$ , and for m-reSGLD is $100\%$.} 
\label{pendulum_confidence_small}
\end{figure}

To test how well the proposed frameworks estimate the uncertainty, we select at random one test trajectory and plot, in Figure~\ref{pendulum_confidence_small}, the estimated confidence intervals. Figure~\ref{pendulum_confidence_small} shows that the proposed frameworks capture the true test trajectory within their confidence interval, \ie $e_3=100\%$. Adam, on the other hand, not only fails to capture the whole true trajectory ($e_3=71\%$), but also overestimates the predictive uncertainty by producing a wider confidence band in some regions.

\begin{figure}[h!]
\centering
\includegraphics[scale = 0.5]{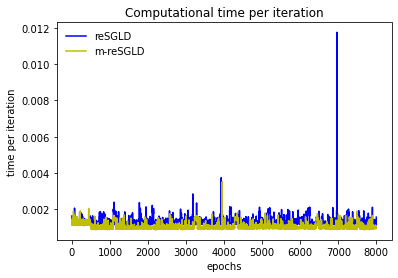}
\caption{Pendulum example with noise $\mathcal{N}(0, 0.01^2)$. Computer running time per iteration for reSGLD and m-reSGLD. The average computation time of m-reSGLD is around $79.06\%$ of reSGLD framework.} 
\label{pendulum_time_small}
\end{figure}

To this end, we compare the computational cost of reSGLD and m-reSGLD, we plot in Figure~\ref{pendulum_time_small} the per-iteration running time for both frameworks. We observe that the per-iteration time for m-reSGLD is smaller than that for reSGLD. Thus, the computational cost for m-reSGLD represents $79.06\%$ of the computational cost for reSGLD. Such a reduction is achieved by m-reSGLD without compromising mean predictive accuracy (see Table~\ref{pendulum_error_small_table}).

\subsubsection{Increased noise}
\begin{figure}[h!]
\centering
\includegraphics[scale = 0.4]{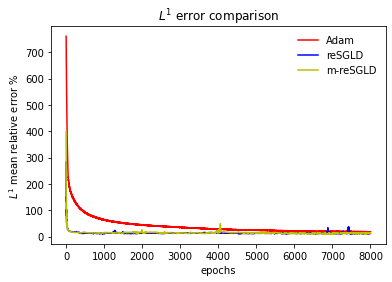}
\includegraphics[scale = 0.4]{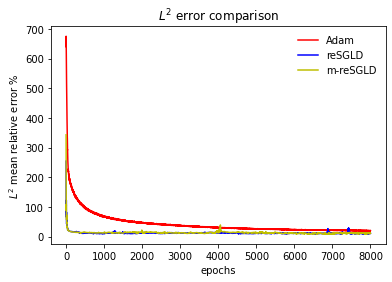}
\includegraphics[scale = 0.4]{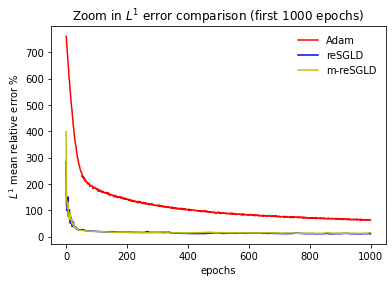}
\includegraphics[scale = 0.4]{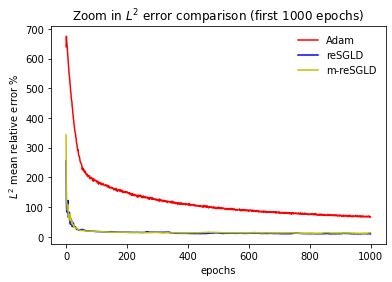}
\includegraphics[scale = 0.4]{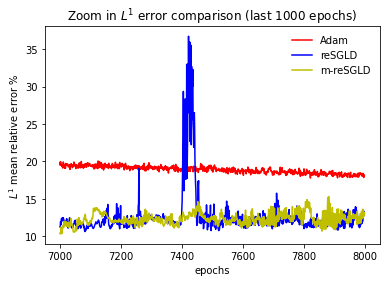}
\includegraphics[scale = 0.4]{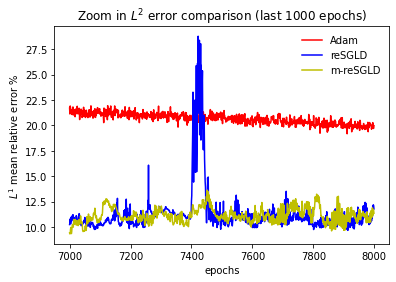}
\caption{Convergence results for the gravity-pendulum with increased noise $\mathcal{N}(0, 0.05^2)$. The $e_1$ errors are shown on the left while $e_2$ errors are shown on the right. Top row: entire history. Middle row: zoom in the first $1,000$ epochs. Last row: zoom in the last $1,000$ epochs. The mean errors after the burn-in are presented in Table~\ref{pendulum_error_big}.} 
\label{pendulum_error_big}
\end{figure}

In this section, we verify the performance of the DeepONets when the variance of noise for the operator targets increases. In particular, we assume the noise follows the normal distribution $\epsilon_i \sim \mathcal{N}(0,0.05^2)$. Figure~\ref{pendulum_error_big} and Table~\ref{pendulum_error_big_table} detail the convergence errors for this increased noise scenario. Similar to the anti-derivative example, the convergence results deteriorate when the noise increases. However, compared with Adam, the proposed frameworks still provide improved training convergence and mean prediction performance. 

\begin{table}[h!]
\centering
\begin{tabular}{||c c c||} 
\hline
frameworks & $e_1$ & $e_2$ \\ [0.5ex] 
\hline
Adam & $18.8935$ & $20.6669$ \\ [0.5ex]
\hline
reSGLD  & $12.6139$ & $11.2885$  \\ [0.5ex]
\hline
m-reSGLD   & $12.2996$ & $11.11955$ \\ [0.5ex]
\hline
\end{tabular}
\caption{Gravity-pendulum example with increased noise $\mathcal{N}(0, 0.05^2)$. The mean errors after burn-in for all three frameworks.}
\label{pendulum_error_big_table}
\end{table}

\begin{figure}[h!]
\centering
\includegraphics[scale = 0.35]{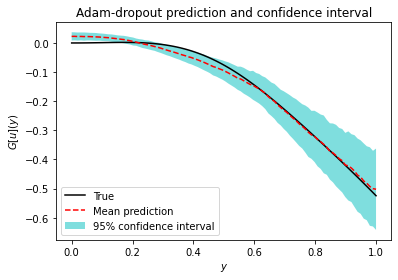}
\includegraphics[scale = 0.35]{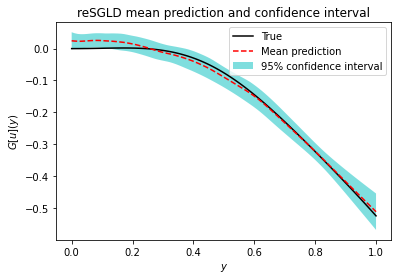}
\includegraphics[scale = 0.35]{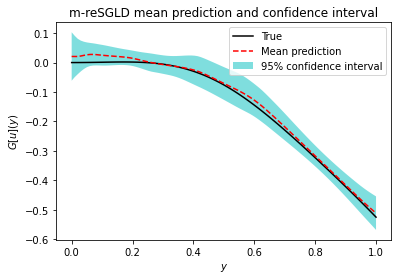}
\caption{Confidence intervals for the gravity pendulum example with noise $\mathcal{N}(0, 1^2)$. Left: Adam. Middle: reSGLD. Right: m-reSGLD.  In this example, the error $e_3$ for Adam is $77\%$, for reSGLD is $93\%$ , and for m-reSGLD is $100\%$.} 
\label{pendulum_confidence_big}
\end{figure}

The constructed confidence intervals for the three frameworks in this increased noise scenario are depicted in Figure~\ref{pendulum_confidence_big}. We observe that the proposed frameworks have a wider confidence band. This is expected since the increase in noise variance leads to larger uncertainty. Adam, however, seems to not change its estimate for this increased noise scenario. Furthermore, compared to Adam $(e_3=77\%)$, we note that the proposed frameworks capture fairly well the true trajectory with their estimated confidence intervals (reSGLD: $e_3=93\%$ and m-reSGLD: $e_3=100\%$). 

\begin{figure}[h!]
\centering
\includegraphics[scale = 0.5]{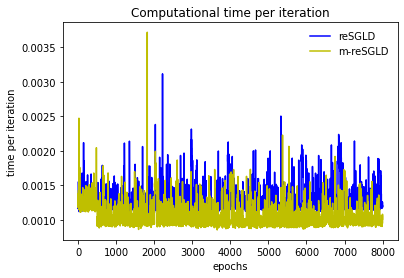}
\caption{Pendulum example with noise $\mathcal{N}(0, 0.05^2)$. Computer running time per iteration for reSGLD and m-reSGLD. The average computation time of m-reSGLD is around $80.8\%$ of reSGLD framework.} 
\label{pendulum_time_big}
\end{figure}

We conclude this gravity pendulum example by illustrating in Figure~\ref{pendulum_time_big} the per-iteration time used by reSGLD and m-reSGLD. The results show that the computational cost of m-reSGLD corresponds to $80.8\%$ of the computational cost of reSGLD. As before, m-reSGLD achieves such a reduction without compromising the mean prediction performance. 

\subsection{Experiment 3}
\begin{figure}[h!]
\centering
\includegraphics[scale = 0.4]{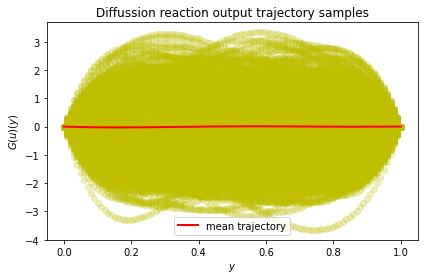}
\caption{Distribution of the operator targets~$G^\dagger(u)(y)$ obtained from the distribution of inputs depicted in Figure~\ref{input_distribution} and used to train the diffusion reaction example.} 
\label{dr_output_distribution}
\end{figure}

In our third experiment, we test the performance of the proposed Bayesian DeepONet with the following diffusion reaction equation,
\begin{subequations}
\begin{align}
    & \frac{\partial s}{\partial t} = D \frac{\partial^2 s}{\partial x^2} + k s^2 + u(x), \quad x \in \Omega,~t \in [0,1],\\
    & u (x, t) = 0, x\in \partial\Omega,\\
    & u(x, 0) = 0, x\in\Omega
\end{align}
\end{subequations}
where $\Omega = [0, 1]$ and $k = -0.01$. Similar to the previous examples, the input samples $u(x)$ are simulated using the GRF (see Figure~\ref{input_distribution}) and discretized using $m=100$ sensors. Figure~\ref{dr_output_distribution} shows the distribution of the true operator targets associated to the input samples~$u$ for the diffusion reaction system.

\subsubsection{Small noise}
\begin{figure}[h!]
\centering
\includegraphics[scale = 0.4]{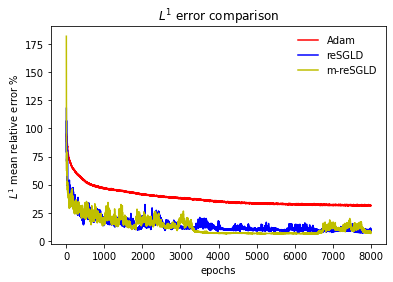}
\includegraphics[scale = 0.4]{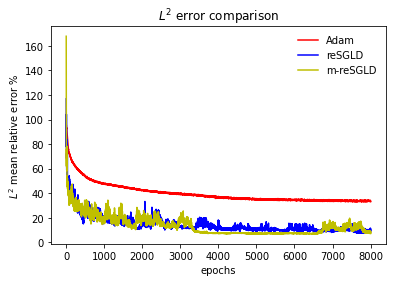}
\includegraphics[scale = 0.4]{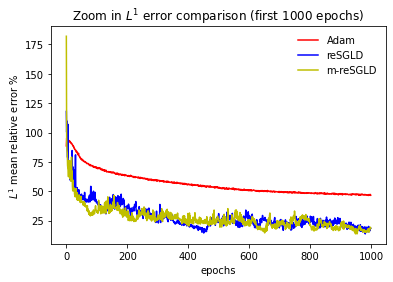}
\includegraphics[scale = 0.4]{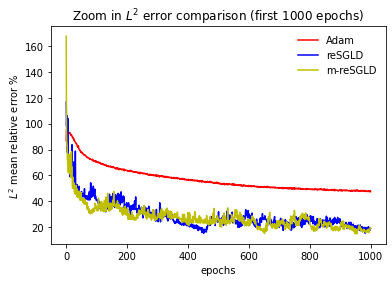}
\includegraphics[scale = 0.4]{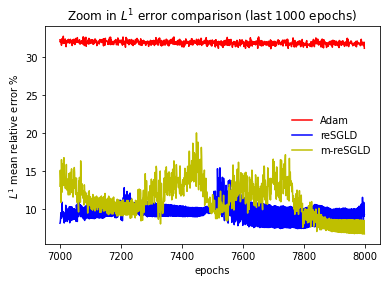}
\includegraphics[scale = 0.4]{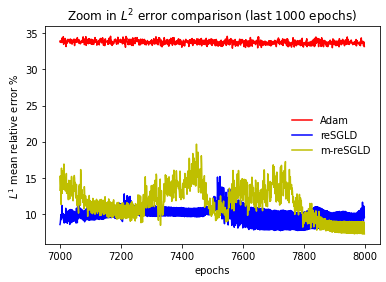}
\caption{Convergence results for the diffusion reaction example with noise $\mathcal{N}(0, 0.01^2)$.  The $e_1$ errors are presented on the left while $e_2$ errors are shown on the right. Top row: entire history. Middle row: zoom in the first $1,000$ epochs. Last row: zoom in the last $1,000$ epochs. The average errors after the burn in are presented in the Table~\ref{dr_error_small_table}.} \label{dr_small_error}
\end{figure}

In this small noise scenario, the noise~$\epsilon_i$ added to the operator targets follows the normal distribution $\mathcal{N}(0, 0.01^2)$. We train the DeepONet for $8,000$ epochs. Figure~\ref{dr_small_error} shows the convergence results of the relative errors~$e_1$ and $e_2$ during training. Also, Table~\ref{dr_error_small_table} presents the average of these relative errors after the burn-in epochs. Compared to the state-of-the-art Adam optimizer, the results show that the proposed frameworks (reSGLD and m-reSGLD) present improved training convergence and mean prediction performance.

\begin{table}[h!]
\centering
\begin{tabular}{||c c c||} 
\hline
frameworks & $e_1$ & $e_2$ \\ [0.5ex] 
\hline
Adam & $31.854993$ & $33.759108$ \\ [0.5ex]
\hline
reSGLD  & $9.644058$ & $ 9.904289$  \\ [0.5ex]
\hline
m-reSGLD   & $11.344359$ & $11.62751$ \\ [0.5ex]
\hline
\end{tabular}
\caption{Diffusion reaction example with noise $\mathcal{N}(0, 0.01^2)$. The mean errors after burn-in for all three frameworks.}
\label{dr_error_small_table}
\end{table}

\begin{figure}[h!]
\centering
\includegraphics[scale = 0.35]{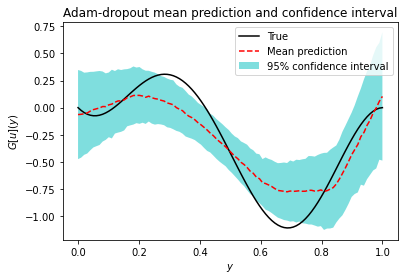}
\includegraphics[scale = 0.35]{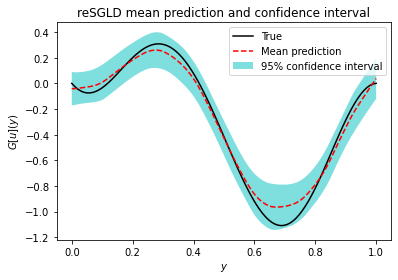}
\includegraphics[scale = 0.35]{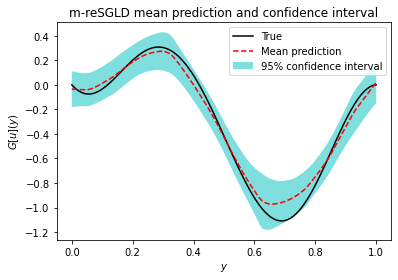}
\caption{Confidence intervals for the diffusion reaction example with noise $\mathcal{N}(0, 0.01^2)$. Left: Adam. Middle: reSGLD. Right: m-reSGLD.  In this example, $e_3$ of Adam is \textcolor{black}{$71\%$, for reSGLD is $100\%$, and for m-reSGLD is $100\%$.}} 
\label{dr_confidence_small}
\end{figure}

We also show in Figure~\ref{dr_confidence_small} the estimated confidence intervals for the three frameworks. The results show that Adam with dropout over-estimates the uncertainty (wider confidence band in some regions), but fails to capture the whole true test trajectory ($e_3 = 71\%$) within the confidence interval. On the other hand, the proposed frameworks, reSGLD and m-reSGLD, provide a more balanced confidence intervals that can capture the whole true trajectory ($e_3 = 100\%$).

\begin{figure}[h!]
\centering
\includegraphics[scale = 0.5]{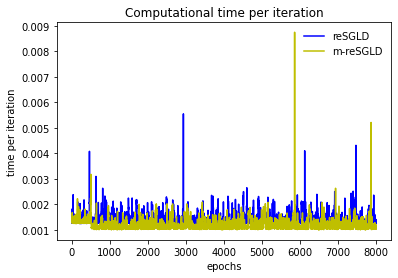}
\caption{Diffusion reaction example with noise $\mathcal{N}(0, 0.01^2)$. Computer running time per iteration for reSGLD and m-reSGLD. The average computation time of m-reSGLD is around \textcolor{black}{$83.0\%$} of reSGLD framework.} 
\label{dr_time_small}
\end{figure}

For this small noise scenario, we plot in Figure~\ref{dr_time_small} the per-iteration time consumed by the proposed frameworks reSGLD and m-reSGLD. Note that by using m-reSGLD we can save $17\%$ of the computational cost of re-SGLD. Furthermore, m-reSGLD achieves such a reduction without sacrificing the mean predictive performance (see Table~\ref{dr_error_small_table}).

\subsubsection{Increased noise}
\begin{figure}[h!]
\centering
\includegraphics[scale = 0.4]{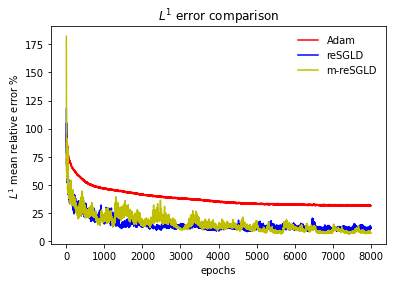}
\includegraphics[scale = 0.4]{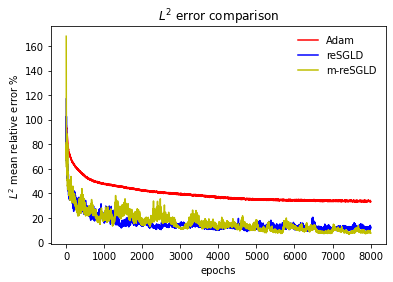}
\includegraphics[scale = 0.4]{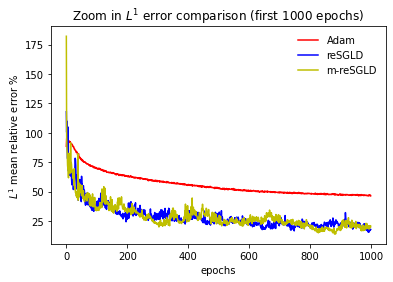}
\includegraphics[scale = 0.4]{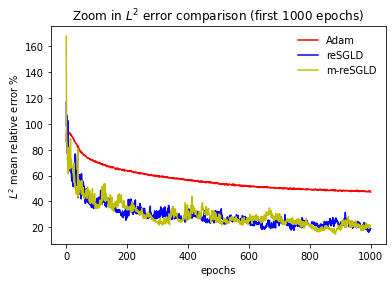}
\includegraphics[scale = 0.4]{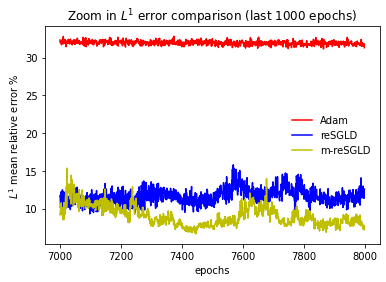}
\includegraphics[scale = 0.4]{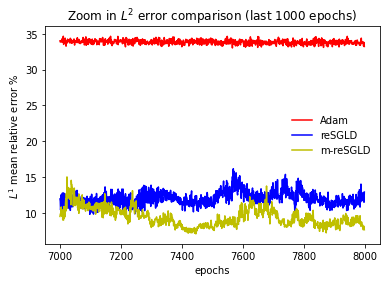}
\caption{Convergence results for the diffusion reaction example with noise $\mathcal{N}(0, 0.1^2)$.  The $e_1$ errors are presented on the left while $e_2$ errors are shown on the right. Top row: entire history. Middle row: zoom in the first $1,000$ epochs. Last row: zoom in the last $1,000$ epochs. The average errors after the burn in are presented in the Table \ref{dr_error_big_table}.} 
\label{dr_big_error}
\end{figure}

In this section, we increase the noise variance so that the noise follows the normal distribution $\mathcal{N}(0, 0.1^2)$. Figure \ref{dr_big_error} shows the convergence results of $e_1$ and $e_2$ relative errors during training. Moreover, Table \ref{dr_error_big_table} presents the average of these errors after the burn-in. Compared to Adam, we observe that the proposed frameworks present improved mean prediction performance and training convergence. 

\begin{table}[h!]
\centering
\begin{tabular}{||c c c||} 
\hline
frameworks & $e_1$ & $e_2$ \\ [0.5ex] 
\hline
Adam & $31.954806$ & $33.850071$ \\ [0.5ex]
\hline
reSGLD  & $11.670668$ & $12.043854$  \\ [0.5ex]
\hline
m-reSGLD   & $9.040757$ & $ 9.404453$ \\ [0.5ex]
\hline
\end{tabular}
\caption{Diffusion reaction example with noise $\mathcal{N}(0, 0.1^2)$. The mean errors after burn-in for all three frameworks.}
\label{dr_error_big_table}
\end{table}

\begin{figure}[h!]
\centering
\includegraphics[scale = 0.35]{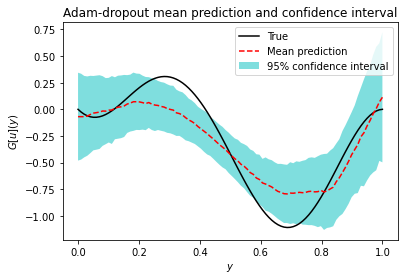}
\includegraphics[scale = 0.35]{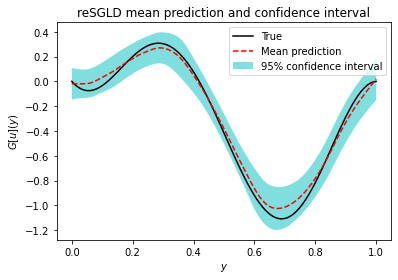}
\includegraphics[scale = 0.35]{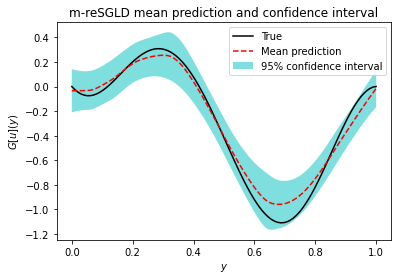}
\caption{Confidence intervals for the diffusion reaction example with noise $\mathcal{N}(0, 0.1^2)$. Left: Adam. Middle: reSGLD. Right: m-reSGLD.  In this example, $e_3$ of Adam is $68\%$, for reSGLD is $100\%$ , and for m-reSGLD is $100\%$.} 
\label{dr_confidence_big}
\end{figure}

The estimated confidence intervals are depicted in Figure~\ref{dr_confidence_big}. Clearly, the proposed frameworks capture much better the true trajectory in their confidence intervals $(e_3=100\%)$. Adam, on the other hand, seems affected by the noise and fails to capture the whole trajectory $(e_3=68\%)$.

\begin{figure}[h!]
\centering
\includegraphics[scale = 0.5]{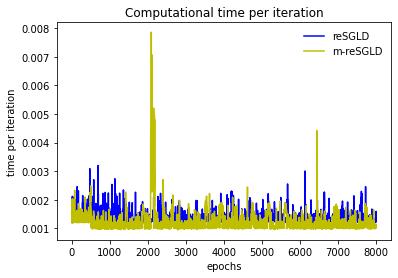}
\caption{Diffusion reaction example with noise $\mathcal{N}(0, 0.1^2)$. Computer running time per iteration for reSGLD and m-reSGLD. The average computation time of m-reSGLD is around \textcolor{black}{$81.5\%$} of reSGLD framework.} 
\label{dr_time_big}
\end{figure}

We conclude this example by plotting the per-iteration time for reSGLD and m-reSGLD (Figure~\ref{dr_time_big}). We conclude that m-reSGLD uses only $81.5\%$ of the computational resources used by reSGLD. Such a reduction is achieved by the proposed m-reSGLD with (remarkably) improved prediction performance \textcolor{black}{(see Table~\ref{dr_error_big_table})}.

\subsection{Experiment 4}
\begin{figure}[h!]
\centering
\includegraphics[scale = 0.5]{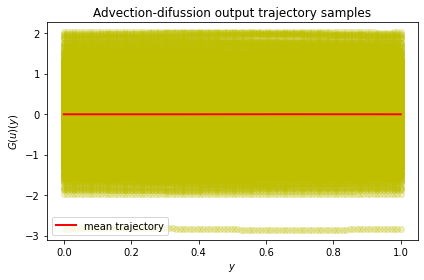}
\caption{Distribution of the operator targets~$G^\dagger(u)(y)$ obtained from the distribution of inputs depicted in Figure~\ref{input_distribution} and used to train the advection-diffusion example.} 
\label{advd_output_distribution}
\end{figure}

In our final experiment, we verify the effectiveness of the proposed Bayesian DeepONet using the following advection-diffusion system with periodic boundary conditions,
\begin{align}
\frac{\partial s}{\partial t} + \frac{\partial s}{\partial x} - D \frac{\partial^2 s}{\partial x^2} &= 0, \qquad x \in \Omega, t \in [0,T],  \\
s(x,0) &= u(x),
\end{align}
where $\Omega = [0,1]$, $T=1$, $D= 0.1$. Here $u(x)$ is an initial condition sampled from a GRF with an RBF kernel of length scale~$l=0.2$ and discretized using $m=100$ sensors. We show in Figure~\ref{advd_output_distribution} the distribution of the true operator targets associated with the input samples~$u$.

\subsubsection{Small noise}
\begin{figure}[h!]
\centering
\includegraphics[scale = 0.4]{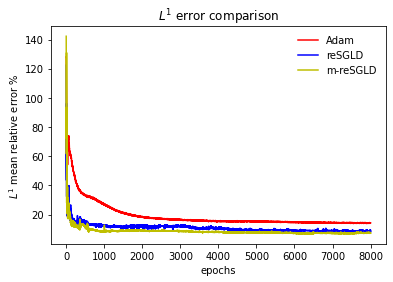}
\includegraphics[scale = 0.4]{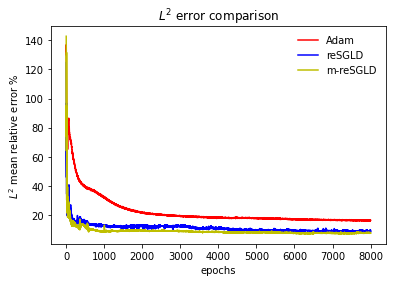}
\includegraphics[scale = 0.4]{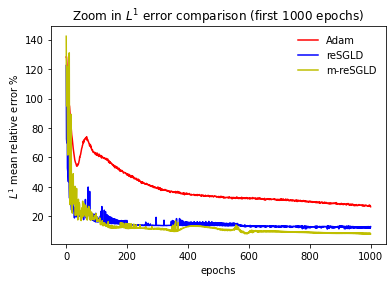}
\includegraphics[scale = 0.4]{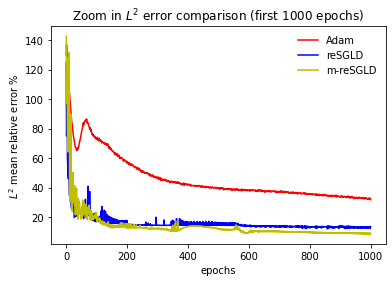}
\includegraphics[scale = 0.4]{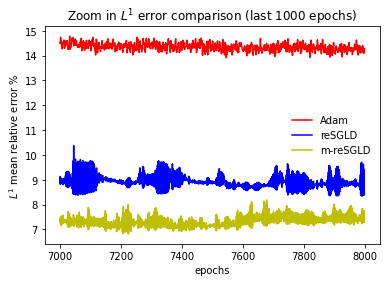}
\includegraphics[scale = 0.4]{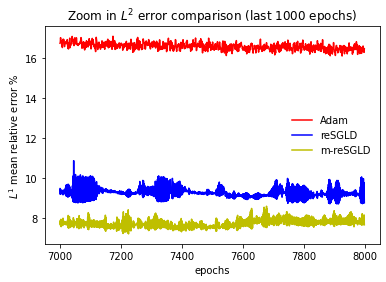}
\caption{Convergence results for the gravity-pendulum with noise $\mathcal{N}(0, 0.01^2)$. The $e_1$ errors are shown on the left while $e_2$ errors are shown on the right. Top row: entire history. Middle row: zoom in the first $1,000$ epochs. Last row: zoom in the last $1,000$ epochs. \textcolor{black}{The average errors after the burn in are presented in the Table~\ref{advd_error_small_table}.}} 
\label{advd_error_small}
\end{figure}

In the small noise scenario, we assume the training operator targets are corrupted by noise that follows the normal distribution~$\epsilon_i \sim \mathcal{N}(0, 0.01^2)$. We train the DeepONets for $8,000$ epochs. Figure \ref{advd_error_small} shows the convergence errors~$e_1$ and $e_2$  and Table \ref{advd_error_small_table} presents the average of these errors after the burn-in epochs. From Figure \ref{advd_error_small} and Table \ref{advd_error_small_table}, we conclude that reSGLD and m-reSGLD have faster training converge and better mean prediction performance than Adam. 

\begin{table}[h!]
\centering
\begin{tabular}{||c c c||} 
\hline
frameworks & $e_1$ & $e_2$ \\ [0.5ex] 
\hline
Adam & $14.3524$ & $16.5843$ \\ [0.5ex]
\hline
reSGLD  & $8.9257$ & $9.2734$  \\ [0.5ex]
\hline
m-reSGLD   & $7.3409$ & $7.7345$ \\ [0.5ex]
\hline
\end{tabular}
\caption{Advection-diffusion example with noise $\mathcal{N}(0, 0.01^2)$. The mean errors after burn-in for all three frameworks.}
\label{advd_error_small_table}
\end{table}

\begin{figure}[h!]
\centering
\includegraphics[scale = 0.35]{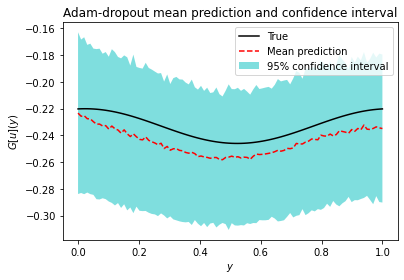}
\includegraphics[scale = 0.35]{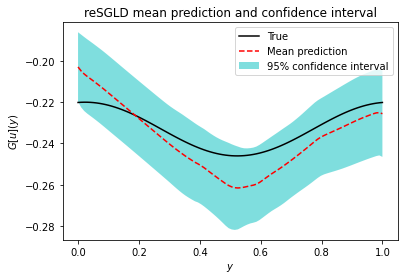}
\includegraphics[scale = 0.35]{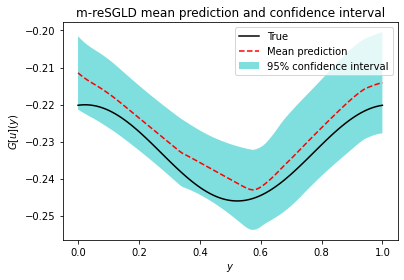}
\caption{Confidence intervals for the advection-diffusion example with noise $\mathcal{N}(0, 0.01^2)$. Left: Adam. Middle: reSGLD. Right: m-reSGLD.  In this example, the error $e_3$ for Adam is $100\%$, for reSGLD is $99\%$ , and for m-reSGLD is $100\%$.} 
\label{advd_confidence_small}
\end{figure}

To illustrate how well our Bayesian DeepONet estimates the uncertainty, we select one test trajectory and plot, in Figure \ref{advd_confidence_small}, the estimated confidence interval by the three frameworks. We observe that the three methods capture the true trajectory $(e_3=100\%)$. Adam, however, seems to over-estimate the uncertainty, yielding a wider confidence band. 

\begin{figure}[h!]
\centering
\includegraphics[scale = 0.5]{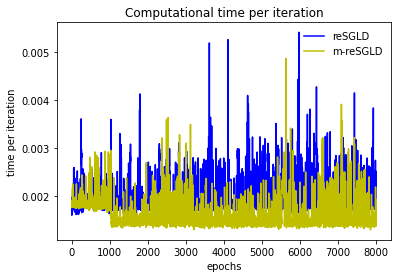}
\caption{\textcolor{black}{Advection-diffusion example with noise $\mathcal{N}(0, 0.01^2)$. Computer  running  time  periteration for m-reSGLD and reSGLD. The average computation time of m-reSGLD is around $79.44\%$ of reSGLD framework. }} 
\label{advd_time_small}
\end{figure}

We conclude this section by plotting (Figure \ref{advd_time_small}) the per-iteration time consumed by the proposed frameworks reSGLD and m-reSGLD. We observe that The average computation time of m-reSGLD is around $79.44\%$ of reSGLD framework. The proposed accelerated training regime m-reSGLD achieves such a reduction even with improved mean prediction performance (see Table~\ref{advd_error_small_table}).

\subsubsection{Increased noise}
\begin{figure}[h!]
\centering
\includegraphics[scale = 0.4]{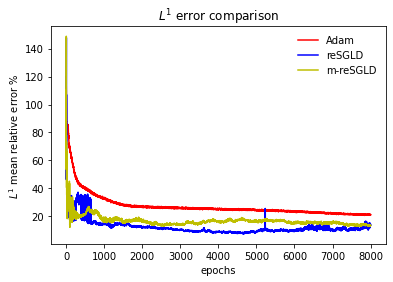}
\includegraphics[scale = 0.4]{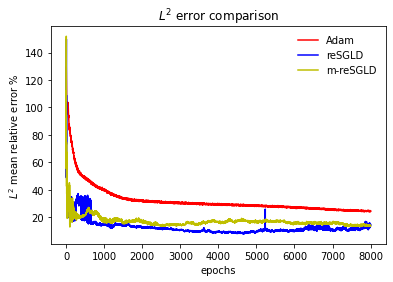}
\includegraphics[scale = 0.4]{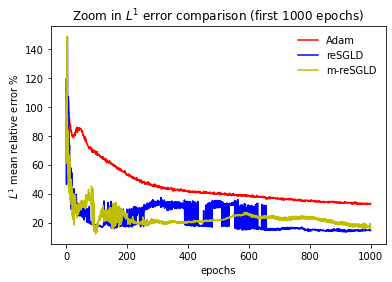}
\includegraphics[scale = 0.4]{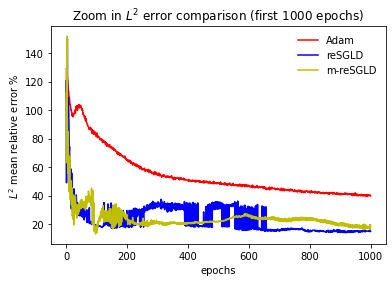}
\includegraphics[scale = 0.4]{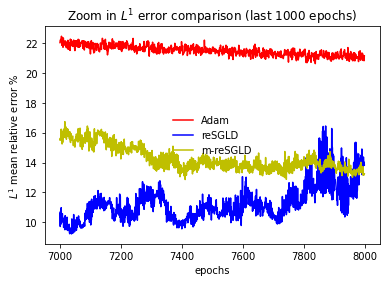}
\includegraphics[scale = 0.4]{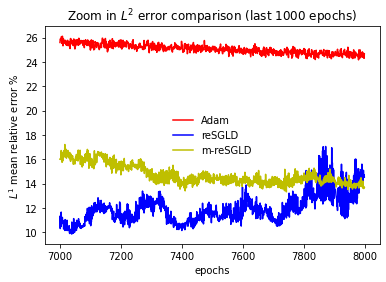}
\caption{Convergence results for the advection-diffusion with increased noise $\mathcal{N}(0, 0.1^2)$. The $e_1$ errors are shown on the left while $e_2$ errors are shown on the right. Top row: entire history. Middle row: zoom in the first $1,000$ epochs. Last row: zoom in the last $1,000$ epochs. \textcolor{black}{The average errors after the burn in are presented in the Table \ref{advd_error_big_table}}.} 
\label{advd_error_big}
\end{figure}

Finally, we verify the performance of the proposed Bayesian DeepONet when the variance of the noise for the operator targets increases. In particular, for the advection-diffusion example, we assume the increased noise follows the normal distribution~$\epsilon_i \sim \mathcal{N}(0,0.1^2)$. Figure \ref{advd_error_big} and Table \ref{advd_error_big_table} detail the convergence errors for this increased noise scenario. Similar to our previous experiments, the proposed frameworks present improved training convergence and better mean prediction performance than Adam. This shows that the proposed frameworks handle noise more effectively.

\begin{table}[h!]
\centering
\begin{tabular}{||c c c||} 
\hline
frameworks & $e_1$ & $e_2$ \\ [0.5ex] 
\hline
Adam & $21.5141$ & $25.0701$\\ [0.5ex]
\hline
reSGLD  & $11.2455$ & $11.9178$  \\ [0.5ex]
\hline
m-reSGLD   & $14.2553$ & $14.7162$ \\ [0.5ex]
\hline
\end{tabular}
\caption{Advection-diffusion example with increased noise $\mathcal{N}(0, 0.1^2)$. The mean errors after burn-in for all three frameworks.}
\label{advd_error_big_table}
\end{table}

\begin{figure}[h!]
\centering
\includegraphics[scale = 0.35]{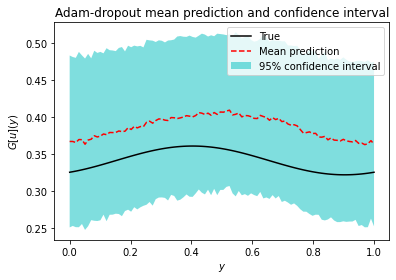}
\includegraphics[scale = 0.35]{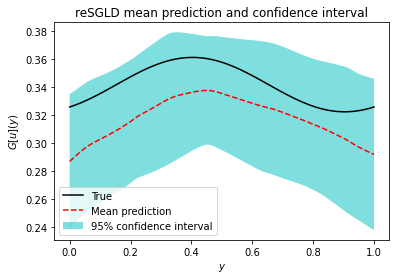}
\includegraphics[scale = 0.35]{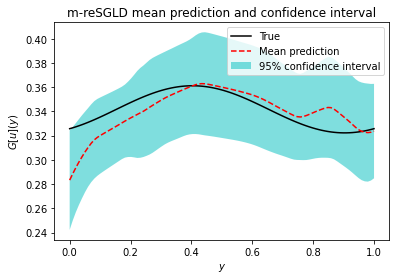}
\caption{Confidence intervals for the advection-diffusion example with noise $\mathcal{N}(0, 1^2)$. Left: Adam. Middle: reSGLD. Right: m-reSGLD.  In this example, the error $e_3$ for Adam is $100\%$, for reSGLD is $100\%$ , and for m-reSGLD is $100\%$.} 
\label{advd_confidence_big}
\end{figure}

We verify next how well the DeepONets capture uncertainty in this increased noise scenario. To this end, we select at random a test trajectory and plot the estimated confidence intervals in Figure \ref{advd_confidence_big}. The results show that, as expected, the proposed frameworks are more uncertain in this increased noise scenario, capturing the whole true test trajectory ($e_3=100\%$). Adam, on the other hand, seems to provide the same estimate of uncertainty as for the small noise scenario. 

We also plot the per-iteration time used by re-SGLD and m-reSGLD in Figure \ref{advd_v90_time}. We observe that m-reSGLD is more efficient; it consumes only $80.2\%$ of the computational cost of re-SGLD. 
\begin{figure}[h!]
\centering
\includegraphics[scale = 0.5]{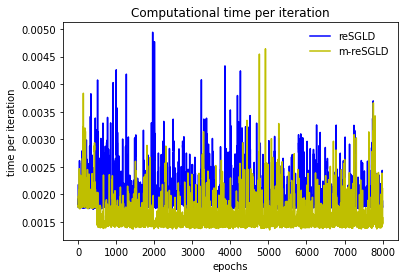}
\caption{Advection-diffusion example with noise $\mathcal{N}(0, 0.1^2)$. Computer running time per iteration for m-reSGLD and reSGLD. The average computation time of m-reSGLD is around $80.2\%$ of reSGLD framework.} 
\label{advd_time_big}
\end{figure}

To conclude this section, we would like to remark that for all the examples above, we can conclude that reSGLD and m-reSGLD converge much faster than Adam during training, with improved prediction accuracy. The proposed frameworks, reSGLD and m-reSGLD, are able to handle the noise in the training data and capture the true solutions within the estimated confidence intervals. The performance of Adam, on the other hand, deteriorates for the increased noise scenario. The benefit of the accelerated training regime of m-reSGLD is that it saves up to 25\% of the computational cost of reSGLD. Such savings are achieved without deteriorating predictive performance. In fact, in some of the examples, m-reSGLD (remarkably) has even better prediction accuracy than reSGLD.

\section{Conclusion} \label{sec:conclusion}
In this work, we study the problem of approximating the solution operator of parametric PDEs using the Deep Operator Network~(DeepONet) framework and noisy data. To this end, we propose using the Bayesian framework of replica-exchange Langevin diffusion. The replica-exchange ability to escape local minima (using one particle to explore and another to exploit the loss function landscape of DeepONets) enables the proposed framework to provide (1) improved training convergence in noisy scenarios, and (2) (remarkably) enhanced predictive capability compared to vanilla DeepONets trained using state-of-the-art gradient-based optimization (\eg Adam). We also demonstrate that the proposed framework effectively estimates the uncertainty of approximating solution trajectories of parametric PDEs. To reduce the computational cost of replica, due to the doubled computational cost of using two particles, we propose an accelerated replica-exchange algorithm that randomly chooses to train one of the DeepONet sub-networks while keeping the other fixed. Using a series of systematically designed experiments, we demonstrate that this accelerated replica-exchange algorithm can save up to 25\% of the original replica-exchange algorithm's original cost without compromising the DeepONet's performance. In our future work, we will study the theoretical aspect of this accelerated replica-exchange framework. 

\bibliographystyle{abbrv}
\bibliography{references}
\end{document}